\documentclass[leqno,12pt]{article} 

\usepackage{amssymb,amsfonts,color}
\usepackage{amsmath,latexsym,mathrsfs}
\usepackage{fancyhdr,verbatim}
\usepackage{graphicx,psfrag}

\newtheorem{theorem}{Theorem}[section]

\newtheorem{lemma}[theorem]{Lemma}

\newtheorem{definition}[theorem]{Definition}
\newtheorem{remark}{Remark}
\newcommand{\ep}{\varepsilon}

\newcommand{\R}{{\mathbb R}}
\newcommand{\N}{{\mathbb N}}

\def\<{\langle}
\def\>{\rangle}
\newcommand{\dif}{{\rm d}}

\renewcommand{\div}{\mbox{\rm div}\;\!}

\newcommand{\fy}{{\varphi}}

\newcommand{\dt}{\partial_{t}}

\newcommand{\cond}{{k}}
\newcommand{\Peclet}{\mbox{\small{\rm Pe}}}

\newcommand{\Stefan}{\mbox{\small{\rm St}}}

\begin{document}

\title{On a phase field model for solid-liquid phase transitions}

\author{S. Benzoni-Gavage\footnote{Universit\'e de Lyon, CNRS UMR 5208 \& Universit\'e Lyon , Institut Camille Jordan, 43 bd du 11 novembre 1918, F-69622 Villeurbanne cedex, France; email: benzoni@math.univ-lyon1.fr} and L. Chupin\footnote{Universit\'e Blaise Pascal \& CNRS UMR 6620, Laboratoire de Math\'ematiques, Campus des C\'ezeaux, B.P. 80026, F-63177 Aubi\`ere cedex, France; email: Laurent.Chupin@math.univ-bpclermont.fr} and D. Jamet\footnote{CEA-Grenoble (DEN/DTP/SMTH), 17, rue des martyrs, F-38054 Grenoble cedex 9, France; email: didier.jamet@cea.fr} and J. Vovelle\footnote{Universit\'e de Lyon, CNRS UMR 5208 \& Universit\'e Lyon , Institut Camille Jordan, 43 bd du 11 novembre 1918, F-69622 Villeurbanne cedex, France; email: vovelle@math.univ-lyon1.fr}}

\date{\today}

\maketitle

\begin{abstract}
A new phase field model is introduced, which can be viewed as a nontrivial generalisation of what is known as the Caginalp model. It involves in particular nonlinear diffusion terms.
By formal asymptotic analysis, it is shown  that in the sharp interface limit it still yields a Stefan-like model with:  1) a generalized Gibbs-Thomson relation telling how much the interface temperature differs from the equilibrium temperature when the interface is moving or/and is curved with surface tension; 2) a jump condition for the heat flux, which turns out to depend on the latent heat and on the velocity of the interface with a new, nonlinear term compared to standard models. From the PDE analysis point of view, the initial-boundary value problem is proved to be locally well-posed in time (for smooth data).
\end{abstract}

\section{Introduction}
Phase field models are widely used in various physical contexts in which a material exhibits two distinct phases.
This is the case for solid-liquid mixtures (e.g.~ice-water or alloys during solidification) or for liquid-vapor mixtures
(e.g.~boiling water), but also for elastic materials subject to martensitic transformations.
The phase field approach is of special interest, in particular  for numerical purposes, when interfaces between the two phases are expected to show complex geometries and topological changes. In phase field models, the `interfaces' are actually viewed as \emph{diffuse interfaces} (see for instance the famous review paper \cite{AndersonMcFaddenWheeler}), {i.e.}~transition regions of nonzero thickness across which a so-called \emph{order parameter} varies smoothly from one to the other of its values in the distinguished phases. Here we are interested in a phase field model designed for solid-liquid mixtures at rest, which consists of an Allen-Cahn type equation for the order parameter coupled with a modified heat equation taking into account both the latent heat and the increase of entropy due to the non-equilibrium situation inside phase-transition regions. This model turns out to be a refined version - in a nontrivial way - of what is known as the Caginalp model \cite{Caginalp}, and it
can also be viewed as a special case of another one designed by Ruyer \cite{Ruyer} for moving liquid-vapor  mixtures.

The aim of this paper is twofold:
1) by formal asymptotic analysis, we show that in the sharp interface limit our model yields a Stefan-like model with
a (generalized) Gibbs-Thomson relation telling how much the interface temperature differs from the equilibrium temperature when the interface is moving or/and is curved with surface tension, together with a jump condition for the heat flux, which turns out to depend on the latent heat and the velocity of the interface with a new, nonlinear term compared to standard models;
2) from the PDE analysis point of view, we prove the local well-posedness of the Cauchy problem and initial-boundary value problems for smooth data. Given that our model displays nonclassical features - it may be seen as a degenerate reaction-diffusion system with nonlinear diffusion - global well-posedness or rough data are not addressed here. 

The mathematical literature on phase-fields equation is extremely vast. In particular, there exist many
extensions of the original Caginalp model developed in \cite{Caginalp}. Let us in particular refer to \cite{GrasselliPetzeltovaSchimperna,BonettiColliFabrizioGilardi09,ColliHilhorstIssardRochSchimperna09,MiranvilleQuintanilla09,CavaterraGalGrasselliMiranville10,GrasselliMiranvilleSchimperna10,MiranvilleQuintanilla10}, which are quite recent papers, and to references therein. We will not attempt to give a minute comparison between 
these references and our work: let us simply emphasize that, up to our knowledge, the model we consider here is {\it distinct} from  
all the models  
considered so far, 
 mainly by the occurrence of the second order quadratic term $\Delta\varphi\partial_t\varphi$ ($\varphi$ being the order parameter) in the equation for the temperature, {\it cf.} the second equation in \eqref{eq:systsimp}. 

The paper is organized as follows. In Section \ref{s:eq} we derive the model and its six-parameters nondimensionalized version.
The sharp interface limit is investigated in Section \ref{s:as}.
Local well-posedness is shown in Section \ref{s:well}.

\section{Phase field equations}\label{s:eq}

\subsection{Derivation and basic properties}\label{ss:der}

\begin{figure}
\centering\includegraphics[height=4cm]{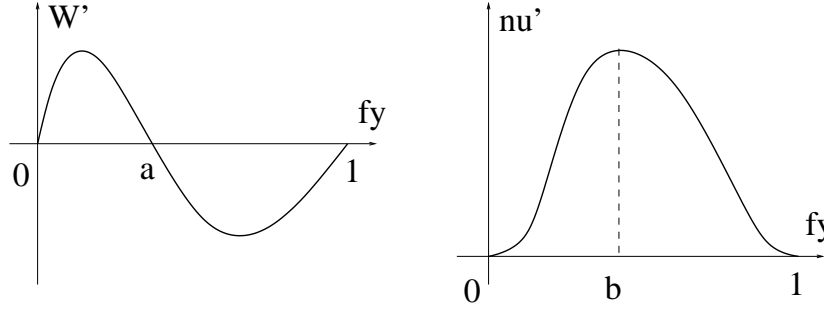}\caption{Derivatives of the double-well potential $W$ and of the entropy coefficient $\nu$.}
\label{fig:Wpnup}
\end{figure}

The model we are going to consider pertains to the so-called second gradient theory. We assume that 
the physical state of a solid-liquid mixture is described by an order parameter $\varphi$ and  its  temperature $T$ in such a way 
that its  free specific energy $f$  depends on $T$, $\varphi$ and also $\nabla \varphi$ in the following way
\begin{equation}\label{eq:freeenergy}
f(T,\varphi,\nabla\varphi)\,=\,\frac{1}{\rho}\,\left(W(\varphi)\,+\,\frac{1}{2}\lambda |\nabla\varphi|^2\right)\,-\,\int_{T_e}^{T} s(\tau,\varphi)\,\dif \tau\,,
\end{equation}
where $\rho$ is the density of the mixture, which will be assumed to be homogeneous and constant,  
$T_e$ is the equilibrium temperature, $\lambda$ is a positive parameter that is supposed to govern the width of solidification/melting fronts,
$W$ is a double-well potential, and $s$ is the specific entropy of the mixture.
More specifically, the order parameter is chosen so that in the pure \emph{phases} 
we have either $\varphi\equiv 1$ (liquid) or 
$\varphi\equiv 0$ (solid), and $W$ is supposed to achieve its global minimum at both $0$ and $1$ and nowhere else.  Furthermore, $s$ is taken to be a convex combination of the entropy in the phases,  depending nonlinearly on the order parameter in the following way
\begin{equation}\label{eq:entropy}
s(T,\varphi)=\nu(\varphi)s_{\mbox{\scriptsize liq}}(T)\,+\,(1-\nu(\varphi))s_{\mbox{\scriptsize sol}}(T)\,,\end{equation}
where $\nu:[0,1]\to [0,1]$ is monotonically increasing.
Typical graphs of the functions $W'$ and $\nu'$ are represented on Fig.~\ref{fig:Wpnup}.
By contrast, in the phase field model of Caginalp \cite{Caginalp},
$\nu$  would be the identity function (hence $\nu'\equiv 1$). 

The \emph{latent heat} of the phase change is by definition
\begin{equation}\label{eq:latentheat}
{\mathcal L}(T)\,:=\,T\,(s_{\mbox{\scriptsize liq}}(T)-s_{\mbox{\scriptsize sol}}(T))\,.
\end{equation}
So another way of writing the free energy is
\begin{multline}\label{eq:freeenergybis}
f(T,\varphi,\nabla\varphi)\\
=\,\frac{1}{\rho}\,\left(W(\varphi)\,+\,\frac{1}{2}\lambda |\nabla\varphi|^2\right)
\,-\,\int_{T_e}^{T} s_{\mbox{\scriptsize sol}}(\tau)\,\dif \tau\,-\,\nu(\varphi)\,\int_{T_e}^{T} \frac{{\mathcal L}(\tau)}{\tau}\,\dif \tau\,,
\end{multline}
and thus the (standard)
\emph{chemical potential} of the mixture is
\begin{equation}\label{eq:chemicalpot}
\mu(T,\varphi):=\frac{\partial f}{\partial \varphi}\,=\,
\frac{1}{\rho}\,W'(\varphi)\,-\,\nu'(\varphi)\,\int_{T_e}^{T}\,\frac{{\mathcal L}(\tau)}{\tau} \,\dif\tau\,.
\end{equation}
Since $W$ has wells at $0$ and $1$, we see that  $\mu\equiv 0$  in both phases whatever the temperature, provided that 
$\nu'$ vanishes at $0$ and $1$ (as on Fig.~\ref{fig:Wpnup}): this would obviously not be the case for a linear $\nu$, as in standard phase field models.
 
The \emph{heat capacity} of the mixture is
\begin{equation}\label{eq:defCp}
C_{p}:=T \frac{\partial s}{\partial T}\,=\,\nu(\varphi)C_{p,\mbox{\scriptsize liq}}(T)\,+\,(1-\nu(\varphi))C_{p,\mbox{\scriptsize sol}}(T).
\end{equation}
 To simplify the analysis, we shall assume that the 
  heat capacities of the liquid $C_{p,\mbox{\scriptsize liq}}=T \partial s_{\mbox{\scriptsize liq}}/\partial T$ and
of the solid $C_{p,\mbox{\scriptsize sol}}=T \partial s_{\mbox{\scriptsize sol}}/\partial T$ have the same constant value $C_0$,
so that $C_{p}$ is also equal to $C_0$. In other words, we shall concentrate on the special case
\begin{equation}\label{eq:simplified}
C_p=C_0\,,\quad {\mathcal L}(T)\,=\,{\mathcal L}_e \frac{T}{T_e}\,,\quad 
s\,=\,s_0+\nu(\varphi)\frac{{\mathcal L}_e}{T_e}\,+\,C_0\,\ln\left(\frac{T}{T_e}\right)\,.\end{equation}

Regardless of that simplifying assumption, we consider the following equations for the evolution of the mixture:
\begin{equation}\label{eq:evolution}
\left\{\begin{array}{l}
\partial_t\varphi\,=\,-\,\kappa\,\mu_g\,,\\  [8pt]
\rho \,C_p\,\partial_t T\,+\,\rho\,\left(\mu_g-T\dfrac{\partial \mu_g}{\partial T}\right)\,\partial_t\varphi\,=\,\div(\cond \nabla T)\,,
\end{array}\right.
\end{equation}
where $\cond> 0$ denotes the \emph{heat conductivity},
$\kappa>0$ denotes the so-called \emph{mobility}, and $\mu_g$ is a
`generalized chemical potential', which merely differs from the standard chemical potential by a second order term:
\begin{equation}\label{eq:defmugen}
\mu_g[T,\varphi]:=\frac{\delta f}{\delta \varphi}\begin{array}[t]{l}\,=\,\mu(T,\varphi)\,-\,\dfrac{1}{\rho}\nabla\cdot (\lambda \nabla \varphi)\,.
\end{array}
\end{equation}
Observe in particular that $\mu_g\equiv 0$ in the phases
($\varphi\equiv 0$ or $\varphi\equiv 1$),
as for $\mu$.
The first equation in \eqref{eq:evolution} is the building block of phase field models,
in which $1/\kappa$ is presumably proportional to a relaxation time for the mixture to return to equilibrium.
Taking into account that
$$C_{p}=T \frac{\partial s}{\partial T}\,,\quad 
 \frac{\partial \mu_g}{\partial T}\,=\, \frac{\partial \mu}{\partial T}\,=\,- \frac{\partial s}{\partial \varphi}\,,$$
the equations in \eqref{eq:evolution} ensure that (for smooth solutions)
\begin{equation}\label{eq:evolutions}
\rho\,\partial_t s\,=\,\div\left(\frac{\cond\nabla T}{T}\right)\, +\,\cond\,\frac{|\nabla T|^2}{T^2}\,+\,\frac{\rho\,\kappa}{T}\;\mu_g^2\,,
\end{equation}
which (formally) means that  the growth of total entropy $\int s(t,x)\,\dif x$ is 
governed by both the conductivity ($\cond$) and the mobility ($\kappa$).

In fact, \eqref{eq:evolution} is specifically designed to have \eqref{eq:evolutions} as well as the (formal) conservation of total energy. More precisely,
the specific energy $e=f+sT$ is  conserved along solutions of \eqref{eq:evolution} in any domain $\Omega\subset \R^d$ such that
\begin{equation}
\label{eq:bcvarphiT}
\int_{\partial \Omega} \lambda \,\nabla \varphi\cdot n\,\partial_t \varphi\,=\,0\;\mbox{and}\; 
\int_{\partial \Omega}  k \,\nabla T\cdot n\,=\,0\,,
\end{equation} where $n$ denotes the normal to $\partial\Omega$.
Indeed, recalling that $s=-\partial f/\partial T$, the first equation in \eqref{eq:bcvarphiT} enables us to write
$$\frac{\dif }{\dif t}\int_\Omega e \,\dif x\,=\,
\int_\Omega \frac{\delta f}{\delta \varphi}\,\partial_t \varphi\,\dif x\,+\,
\int_\Omega T\,\partial_t s\,\dif x\,,
$$
where, by definition of $\mu_g$ and by the first equation in 
\eqref{eq:evolution}, the first integral equals $-\int \kappa \,\mu_g^2\,\dif x$,
which obviously cancels out with the integral coming from the last term in  Eq.~\eqref{eq:evolutions}.
To conclude that $\int_\Omega e(t,x) \,\dif x$ is constant,
 we observe that by  the condition on  $T$ in \eqref{eq:bcvarphiT},
$$\int_\Omega \left(T\,\div\left(\frac{\cond\nabla T}{T}\right)\, +\,\cond\,\frac{|\nabla T|^2}{T}\right)\,\dif x\,=\,0\,.$$

To finish with these general observations, we point out that the equalities in \eqref{eq:bcvarphiT} are easily achieved by means of standard boundary conditions. Namely, the second equality will be implied by a homogeneous Neumann boundary condition on $T$, \emph{i.e.}
$\nabla T\cdot n=0$ (meaning zero heat flux at the boundary: incidentally, one may note that for a nonzero heat flux the total energy will either decrease, due to cooling, or increase, due to heating), and
 either a homogeneous Neumann boundary condition  or 
a Dirichlet condition $\varphi_{|\partial \Omega}\equiv 0$ or $1$ (both values implying $\mu_g\equiv 0$, as already noticed)
will ensure the first one. The appropriate choice of a boundary condition for $\varphi$ is related to the moving contact line problem, which we shall not discuss here.

\subsection{Nondimensionalization}\label{ss:nondim}

The total number of independent physical units used to describe all dependent variables and independent variables in \eqref{eq:evolution}--\eqref{eq:defmugen}  is 
ten
(those of $x$, $t$, $\kappa$, $W$, $\rho$, $T$, ${\mathcal L}$, $\lambda$, $C_p$, $\cond$,
and of course $\varphi$ and $\nu$ do not count because they are already nondimensional),
and the number of fundamental physical units is four ($\mbox{\rm kg}$, $\mbox{\rm m}$, $\mbox{\rm s}$, $\mbox{\rm K}$).
So by elementary dimensional analysis (Buckingham $\pi$ theorem), a nondimensionalized version of 
 \eqref{eq:chemicalpot}--\eqref{eq:evolution}--\eqref{eq:defmugen} requires $10-4=6$ nondimensional parameters.
Below is a possible choice for these parameters, expressed in terms of 
 \begin{itemize}
 \item the density $\rho$,
 \item the equilibrium temperature $T_e$ together with a characteristic temperature difference $\delta\!\,T$,
 \item a length scale $L$,
 \item a characteristic interface thickness $h$,
 \item a time scale $t_0$,
 \item the surface tension $\sigma$,
 \item the latent heat ${\mathcal L}_e$ at $T_e$,
 \item a reference heat capacity $C_0$,
 \item a reference mobility coefficient $\kappa_0$,
 \item a reference heat conductivity $\cond_0$.
  \end{itemize}
  Introducing the  parameters 
  $$\varepsilon:=\frac{h}{L}\,,\;
  \Peclet:= \frac{\rho C_0 L^2}{\cond_0 t_0}\,,\;
  \alpha:=\frac{\kappa_0t_0\sigma}{\rho h}\,,\;
  \theta:=\frac{T_e}{\delta\!\,T}\,,\;
  \beta:=\frac{\sigma}{\rho C_0h\delta\!\,T}\,,\;
  \Stefan:=\frac{C_0\delta\!\,T}{{\mathcal L}_e}\,,
  $$
together with the rescaled variables
  $$\widetilde{x}:=\frac{x}{L}\,,\;\widetilde{t}:=\frac{t}{t_0}\,,\;\widetilde{\kappa}:=\frac{\kappa}{\kappa_0}\,,\;\widetilde{W}:=\frac{W}{\sigma/h}\,,\;\widetilde{T}:=\frac{T-T_e}{\delta\!\,T}\,,\;\widetilde{\mathcal L}:=\frac{{\mathcal L}}{{\mathcal L}_e}\,,\;\widetilde{\lambda}:=\frac{\lambda}{\sigma h}\,,$$
  $$\;\widetilde{C}_p:=\frac{C_p}{C_0}\,,\;\widetilde{\cond}:=\frac{\cond}{\cond_0}\,,$$
we may rewrite \eqref{eq:evolution}--\eqref{eq:defmugen}
as
\begin{equation}\label{eq:evolutionnondim}
\left\{\begin{array}{l}
\partial_{\widetilde{t}}\varphi\,=\,-\,\alpha\,\widetilde{\kappa}\,\widetilde{\mu}_g\,,\\  [8pt]
\widetilde{C}_p\,\partial_{\widetilde{t}} \widetilde{T}\,+\,\beta\,\left(\widetilde{\mu}_g-(\theta+\widetilde{T})\dfrac{\partial \widetilde{\mu}_g}{\partial \widetilde{T}}\right)\,\partial_{\widetilde{t}}\varphi\,=\,\dfrac{1}{\Peclet}\,\nabla_{\widetilde{x}}\cdot (\widetilde{\cond} \nabla_{\widetilde{x}} \widetilde{T})\,,
\end{array}\right.
\end{equation}
\begin{equation}\label{eq:defmugennondim}
\widetilde{\mu}_g[\widetilde{T},\varphi]\,=\,
\widetilde{W}'(\varphi)\,-\,\frac{1}{\beta\Stefan}\,\nu'(\varphi)\,\int_{0}^{\widetilde{T}}\,\frac{\widetilde{{\mathcal L}}(\tau)}{\theta+\tau} \,\dif\tau \,-\,\varepsilon^2\,\nabla_{\widetilde{x}}\cdot (\widetilde{\lambda} \nabla_{\widetilde{x}} \varphi)\,,
\end{equation}
If ${\kappa}$ is supposed to be constant, or similarly if ${\cond}$ is constant, 
we may assume without loss of generality that
$\widetilde{\kappa}\equiv 1$, respectively $\widetilde{\cond}\equiv 1$, in  \eqref{eq:evolutionnondim}. 
The case of $\widetilde{\lambda}$ in \eqref{eq:defmugennondim} is more subtle because 
it depends on the parameters $h$ and ${\lambda}$, the former being arbitrary and the latter not
being accessible to physical measurements. 
Nevertheless, we choose to set $\widetilde{\lambda}\equiv 1$.
In addition, under the simplifying assumptions in \eqref{eq:simplified},  we have in the rescaled variables 
$$\widetilde{C}_p\equiv 1\,,\;\widetilde{{\mathcal L}}\,=\,1+\frac{\widetilde{T}}{\theta}\,,$$
so that in this case the nondimensionalized version \eqref{eq:evolutionnondim}--\eqref{eq:defmugennondim}
 of  \eqref{eq:evolution}--\eqref{eq:defmugen} reads, dropping the tildes for simplicity,
 \begin{equation}\label{eq:evolutionnondimsimp}
\left\{\begin{array}{l}
\partial_t\varphi\,=\,-\,\alpha\,\mu_g\,,\\  [8pt]
\partial_t T\,+\,\beta\,\left(\mu_g-(\theta+T)\dfrac{\partial \mu_g}{\partial T}\right)\,\partial_t\varphi\,=\,\dfrac{1}{\Peclet}\,\Delta T\,,
\end{array}\right.
\end{equation}
\begin{equation}\label{eq:defmugennondimsimp}
\mu_g[T,\varphi]\,=\,
W'(\varphi)\,-\,\frac{1}{\beta\Stefan}\,\nu'(\varphi)\,\frac{T}{\theta} \,-\,\varepsilon^2\,\Delta \varphi\,.
\end{equation}
Plugging \eqref{eq:defmugennondimsimp} into \eqref{eq:evolutionnondimsimp} we get the system
 \begin{equation}\label{eq:systsimp}
\left\{\begin{array}{l}
\partial_t\varphi\,=\,-\alpha\,W'(\varphi)\,+\,\dfrac{\alpha}{\beta\Stefan}\,\nu'(\varphi)\,\dfrac{T}{\theta} \,+\,\alpha\,\varepsilon^2\,\Delta \varphi\,,\\  [8pt]
\partial_t T\,+\,\left(\beta\,W'(\varphi)+\dfrac{1}{\Stefan}\,\nu'(\varphi)\,-\,\beta\,\varepsilon^2\,\Delta \varphi \right)\,\partial_t\varphi\,=\,\dfrac{1}{\Peclet}\,\Delta T\,.
\end{array}\right.
\end{equation}
This resembles the system considered by Caginalp in his seminal paper \cite{Caginalp}, except for two important differences. The first one is that the coefficient of $T$ depends on $\varphi$ in the first equation.
The other one lies in the complicated, second order and nonlinear coefficient of $\partial_t\varphi$ in the second equation, which is -- up to the authors knowledge --, always supposed to be a constant (latent heat) in Caginalp-like models.

For completeness, let us now derive the nondimensional versions of the entropy equation \eqref{eq:evolutions} and of the local conservation law for the energy.
Redefining $s$ as the  nondimensional entropy $s/C_0$, we have from \eqref{eq:simplified} that $$s= \dfrac{1}{\Stefan\theta}\,\nu(\varphi)\,+\,\ln(T+\theta)$$
up to a harmless additive constant.
Then the nondimensionalized version of the entropy equation \eqref{eq:evolutions} is
\begin{equation}\label{eq:evolutionsnondim}
\partial_t s\,=\,\frac{1}{\Peclet}\,\div\left(\frac{\nabla T}{T+\theta}\right)\, +\,\frac{1}{\Peclet}\,\,\frac{|\nabla T|^2}{(T+\theta)^2}\,+\,\frac{\beta}{\alpha(T+\theta)}\,(\partial_t\varphi)^2\,.
\end{equation}
Regarding the nondimensionalized energy 
$$e
:=f
+\,(T+\theta)\,s
\,=\,
T\,+\,\beta \,W(\varphi)\,+\,\dfrac{1}{\Stefan}\,\nu(\varphi)\,+\,\frac{1}{2}\,\beta\,\varepsilon^2 \, |\nabla \varphi|^2$$
we easily find  the conservation law
\begin{equation}\label{eq:evolutionenondim}
\partial_t e\,=\,\frac{1}{\Peclet}\,\Delta T\,+\,\beta\,\varepsilon^2\,\nabla\cdot( (\partial_t \varphi)\,\nabla \varphi)\,.
\end{equation}

\section{Sharp interface limit} \label{s:as}
Our aim here to derive at least formally a physically realistic, asymptotic limit of the system  \eqref{eq:systsimp} when 
the width of interfaces tends to zero, either because of a physical scaling or for other reasons related to the actual values of the six nondimensional parameters $\alpha$, $\beta$, $\varepsilon$, $\theta$, $\Peclet$, and $\Stefan$.
More precisely, we are going to show in what follows that for  suitable relationships between those parameters,
the system  \eqref{eq:systsimp} formally tends to the Stefan-like model \eqref{eq:genStefanphys} (see p.~\pageref{eq:genStefanphys} hereafter) when $\varepsilon$ goes to zero. (Recall that $\varepsilon$ is the parameter governing the typical width of interfaces.) Before entering into details, let us emphasize that the sharp interface model  in \eqref{eq:genStefanphys} naturally involves the heat equation in the phases, and two sorts of conditions at interfaces, namely
\begin{itemize}
\item
a (generalized) Gibbs-Thomson relation giving the interface temperature in terms of the surface tension, the mean curvature and the velocity of the interface,
\item a jump condition for the heat flux across the interface, in terms of the latent heat and of the velocity of the interface,
the later dependence being nonlinear (quadratic).
\end{itemize}
This should be of interest to discuss the physical validity of  \eqref{eq:systsimp}.
\subsection{Formal asymptotics}\label{ss:deras}
\begin{figure}
\centering\includegraphics[width=7cm]{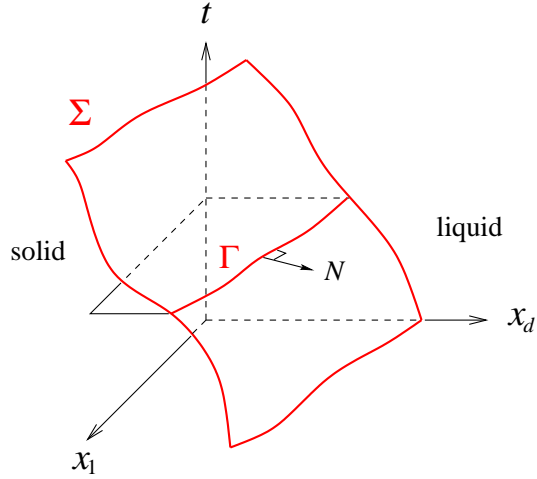}
\caption{Sharp interface configuration}
\label{fig:sharp}
\end{figure}

For convenience, we rewrite \eqref{eq:systsimp} as
 \begin{equation}\label{eq:systsimpbis}
\left\{\begin{array}{l}
\hat{\alpha}\,\partial_t\varphi\,=\,\varepsilon^2\,\Delta \varphi\,-\,W'(\varphi)\,+\,\gamma\,\nu'(\varphi)\,{T}\,,\\  [8pt]
\hat{\beta}\,\partial_t T\,=\,\delta\,\Delta T\,-\,\gamma\,(T+\theta)\,\partial_t\nu(\varphi)\,+\,\hat\alpha\,(\partial_t\varphi)^2\,,
\end{array}\right.
\end{equation}
with
\begin{equation}\label{eq:newparam}
\hat{\alpha}:=\frac{1}{\alpha}\,,\quad \hat{\beta}:=\frac{1}{\beta}\,,\quad \gamma:=\dfrac{1}{\beta\Stefan\theta}\,,\quad 
\delta:= \dfrac{1}{\beta\Peclet}\,.
\end{equation}
The six nondimensional parameters in \eqref{eq:systsimpbis}
are now
$\hat\alpha$, $\hat\beta$, $\gamma$, $\delta$, $\varepsilon$, and $\theta$.
If we go back to the original definitions of $\alpha$, $\beta$, $\Stefan$, $\theta$, and $\Peclet$, we see from  \eqref{eq:newparam} that 
$\hat\alpha$, $\hat\beta$, $\gamma$, $\delta$ are all proportional to the ratio $h/\sigma$ of the interface width and the surface tension, and each of them has its own a physical significance according to the following relationships
$$\begin{array}{ll}\hat\alpha \propto 1/\kappa_0 & \mbox{(relaxation time)}\,,\\
\hat\beta \propto C_0 & \mbox{(heat capacity)}\,,\\
\gamma \propto {\mathcal L}_e & \mbox{(latent heat)}\,,\\
\delta \propto k_0 & \mbox{(thermal conductivity)}\,.\end{array}$$
As regards the sharp interface limit $\varepsilon=h/L \to 0$ at  fixed surface tension $\sigma$, by the observation above it is rather natural to let the four parameters $\hat\alpha$, $\hat\beta$, $\gamma$, and $\delta$ go to zero at least like $\varepsilon$.  If in addition we let the relaxation time go to zero like $\varepsilon$,
we are led to consider 
$$\overline{\alpha}:=\hat\alpha/\varepsilon^2\,,\;\overline{\beta}:=\hat\beta/\varepsilon\,,\;
\overline{\gamma}:=\gamma/\varepsilon\,,\;\overline{\delta}:=\delta/\varepsilon$$
as being fixed. With these definitions, \eqref{eq:systsimpbis} becomes
\begin{equation}\label{eq:systsimpter}
\left\{\begin{array}{l}
\overline{\alpha}\,\varepsilon^2\,\partial_t\varphi\,=\,\varepsilon^2\,\Delta \varphi\,-\,W'(\varphi)\,+\,\overline{\gamma}\,\varepsilon\,\nu'(\varphi)\,{T}\,,\\  [8pt]
\overline{\beta}\,\partial_t T\,=\,\overline{\delta}\,\Delta T\,-\,\overline{\gamma}\,(T+\theta)\,\partial_t\nu(\varphi)\,+\,\overline{\alpha}\,\varepsilon\,(\partial_t\varphi)^2\,,
\end{array}\right.
\end{equation}
The formal limit of the first equation in \eqref{eq:systsimpter} as $\varepsilon\to 0$ gives $W'(\varphi)=0$, which imposes that $\varphi$ takes only the values $0$ (solid phase), $1$ (liquid phase), or $a$ (`metastable' state), while the formal limit of the second equation is
\begin{equation}\label{eq:Tlim}
\overline{\beta}\,\partial_t T\,=\,\overline{\delta}\,\Delta T\,-\,\overline{\gamma}\,(T+\theta)\,\partial_t\nu(\varphi)\,.
\end{equation}
Assume that $T$ is a continuous solution of  \eqref{eq:Tlim}, in which $\varphi$  represents a sharp interface, that is, $\varphi$ is constant and equal to $0$ or $1$ on either side a smooth, moving surface $\Gamma(t)$ as on Fig.~\ref{fig:sharp}. Then by integration by parts in the neighborhood of any point
$(x,t)\in\Sigma:=\{(x,t)\in \R^d \times \R\,,\;x\in \Gamma(t)\}$  we find 
that the gradient of $T$ experiences a discontinuity across $\Gamma(t)$ according to the following relation
\begin{equation}\label{eq:nablaT}
- \overline\gamma (T+\theta) \,v\,=\,\overline\delta \,\left[\nabla T\cdot N\right]
\end{equation}
where $N$ denotes the unit normal to $\Gamma$ pointing to the liquid phase ($\varphi\equiv 1$) and $v$ denotes the speed of $\Gamma$ in the direction $N$.
A linear relation between $v$ and $\left[\nabla T\cdot N\right]$ as in
Equation \eqref{eq:nablaT} is a classical building block in Stefan models for sharp interfaces, see for instance Fig.~1 in \cite{CaginalpChen}.

Of course the formal limits above are not valid in regions where $(\varphi,T)$ experience large variations, and to describe exact 
solutions of \eqref{eq:systsimpter} we need internal layers for diffuse interfaces.
In what follows we adopt the same, multiscale approach as in \cite{CaginalpChen}, where the sharp interface limit was obtained for
the usual Caginalp model. Consider a (smooth) solution $(\varphi^\varepsilon,T^\varepsilon)$ of \eqref{eq:systsimpter} and let $\Gamma^\varepsilon(t)$ be the level surface $\{x\in\R^d\,;\;\varphi^\varepsilon(x,t)=b\}$ (the $b$ where $\nu'$ attains its maximum, which is supposed to best describe the location of the `interface'). We assume that $\Gamma^\varepsilon(t)$ is  smooth, not self-intersecting, and
depends smoothly on $t$ and $\varepsilon$ in  such a way that the signed distance $d^\varepsilon(x,t)$  of $x$ to   $\Gamma^\varepsilon(t)$ is well-defined for $t\in [0,T]$, $\varepsilon\in [0,\varepsilon_0]$, and $x$ in some neighborhood ${\mathscr V}(t)$ of 
$\cup_{\varepsilon \in [0,\varepsilon_0]}\Gamma^\varepsilon(t)$.
By definition of $d^\varepsilon$,
$$N^\varepsilon(x,t):=\nabla d^\varepsilon(x,t)$$ 
is a unit normal vector to $\Gamma^\varepsilon(t)$,
$$v^\varepsilon(x,t):=-\partial_td^\varepsilon(x,t)$$
 is the normal speed of $\Gamma^\varepsilon(t)$ in the direction $N^\varepsilon$, 
and $$H^\varepsilon(x,t):=-\Delta d^\varepsilon(x,t)$$ is  the sum of principal curvatures of  $\Gamma^\varepsilon(t)$ at $x$.
Assume that the functions $d^\varepsilon$, $\varphi^\varepsilon$, and $T^\varepsilon$ admit asymptotic expansions of the form
$$d^\varepsilon(x,t) \sim \sum_{i=0}^{\infty} \varepsilon^i\,d_i(x,t)\,,$$
$$\varphi^\varepsilon(x,t) \sim \sum_{i=0}^{\infty} \varepsilon^i\,\varphi_i(x,t,d^\varepsilon(x,t)/\varepsilon)\,,$$
$$T^\varepsilon(x,t) \sim \sum_{i=0}^{\infty} \varepsilon^i\,T_i(x,t,d^\varepsilon(x,t)/\varepsilon)\,,$$
as $\varepsilon\to 0$, for $x\in {\mathscr V}(t)$.
We recall that the notation 
$d^\varepsilon \sim \sum_{i=0}^{\infty} \varepsilon^i\,d_i$ means that
for all $I\in \N$, 
$d^\varepsilon - \sum_{i=0}^{I} \varepsilon^i\,d_i\,=\,o(\varepsilon^{I+1})$.
We shall denote by $z=d^\varepsilon(x,t)/\varepsilon$ the rescaled variable in the normal direction to $\Gamma^\varepsilon(t)$.
For any smooth function $(x,t,z)\mapsto \widetilde{F}(x,t,z)$, the derivatives of the function
$F^\varepsilon:(x,t)\mapsto  \widetilde{F}(x,t,d^\varepsilon(x,t)/\varepsilon)$ are given by
$$\partial_t F^\varepsilon= \partial_t \widetilde{F} +{\varepsilon^{-1}}\,(\partial_z \widetilde{F})\,\partial_t d^\varepsilon\,,$$
$$\nabla  F^\varepsilon= \nabla   \widetilde{F} +{\varepsilon^{-1}}\,(\partial_z \widetilde{F})\,\nabla  d^\varepsilon\,,$$
$$\Delta   F^\varepsilon= \Delta    \widetilde{F} +{\varepsilon^{-1}}\,\left( (\partial_z \widetilde{F})\,\Delta   d^\varepsilon\,+\,2\,\nabla  d^\varepsilon \cdot \nabla \partial_z \widetilde{F}\right)\,+\,\varepsilon^{-2}\,\partial_z^2\widetilde{F}\,,$$
where the differential operators $\nabla$ and $\Delta$ concern only the  variable $x$.
Hence the system \eqref{eq:systsimpter} for $(\varphi,T)=(\varphi^\varepsilon,T^\varepsilon)$ is equivalent to the following one, 
evaluated at $z=d^\varepsilon(x,t)/\varepsilon$:
\begin{equation}\label{eq:systext}
\left\{\begin{array}{l}
\partial_z^2\widetilde\varphi\,-\,W'(\widetilde\varphi)\,=\,\begin{array}[t]{l}\varepsilon\,\left((\overline\alpha\partial_td^\varepsilon-\Delta d^\varepsilon)\partial_z\widetilde\varphi-2\nabla  d^\varepsilon \cdot \nabla \partial_z \widetilde{\varphi}-\overline\gamma\,\nu'(\varphi)\,\widetilde{T}\right)\\ [10pt] \,+\,\varepsilon^2\,(\overline\alpha \partial_t\widetilde\varphi\,-\,\Delta\widetilde\varphi)\,,\end{array}\\ [30pt]
\overline\delta\,\partial_z^2\widetilde{T}\,=\,\begin{array}[t]{l}\varepsilon\,\left(
(\overline\beta\partial_td^\varepsilon-\overline\delta\Delta d^\varepsilon)\partial_z\widetilde{T}\,-\,
2\overline\delta\nabla d^\varepsilon\cdot \nabla\partial_z\widetilde{T}\right.\\ [10pt]
\left.\quad +\,\overline\gamma(\widetilde{T}+\theta)\nu'(\widetilde\varphi)(\partial_td^\varepsilon)\partial_z\widetilde\varphi\,-\,2\,
\overline\alpha\,(\partial_td^\varepsilon)^2(\partial_z\widetilde\varphi)^2\right)\\ [10pt]
+\,\varepsilon^2\,\left(\overline\beta\partial_t\widetilde{T}-\overline\delta\widetilde{T}+\overline\gamma(\widetilde{T}+\theta)\partial_t\nu(\widetilde\varphi)\,-\,\overline\alpha (\partial_td^\varepsilon)(\partial_t\widetilde\varphi)\partial_z\widetilde\varphi\right)\\ [10pt]
-\,\varepsilon^3\,\overline\alpha\,(\partial_t\widetilde\varphi)^2
\,.\end{array}
\end{array}\right.
\end{equation}
We expect that $\Gamma^\varepsilon(t)$ converges to 
$\Gamma^0(t)$, the level set $\{x\in\R^d\,;\;d_0(x,t)=0\}$.
Off $\Gamma^0(t)$, $d^\varepsilon(\cdot,t)/\varepsilon$ tends to $\pm\infty$ as $\varepsilon\to 0$, so that we shall need extensions
of $(\widetilde\varphi,\widetilde{T})$ for all values of $z\in (-\infty,+\infty)$. However, the only constraint is that 
 \eqref{eq:systext} holds at $z=d^\varepsilon(x,t)/\varepsilon$, which means that we can add to the equations any `reasonable' function of $(x,t,d^\varepsilon(x,t)-\varepsilon z)$ that vanishes when its last variable equals zero.
 This observation will be used in a crucial way to deal with the equation on $\widetilde{T}$. 

Retaining only the $\varepsilon^0$ terms in the first equation of \eqref{eq:systext} we get 
\begin{equation}\label{eq:phi0}\partial_z^2\varphi_0-W'(\varphi_0)\,=\,0\,,\end{equation}
which is the standard equation for a stationary diffuse interface connecting 
$0$ at $z=-\infty$ to $1$ at $z=+\infty$ (or vice-versa). A straightforward phase portrait analysis shows that there is a 
unique such $\varphi_0$ satisfying $\varphi_0(0)=b$. In particular, there is no degree of freedom for $\varphi_0$ to depend on $(x,t)$.

To the next order, using that $\varphi_0$ is independent of $x$, we obtain from the factors of $\varepsilon^1$ the equation
\begin{equation}\label{eq:phi1}
(\partial_z^2-W''(\varphi_0))\,\varphi_1\,=\,\,(\overline\alpha\,\partial_t d_0\,-\,\Delta d_0)\,\partial_z \varphi_0\,-\,\overline\gamma \nu'(\varphi_0)\,T_0\,.
\end{equation}
Since $\varphi_0$ tends to its endstates exponentially fast, the right-hand side of \eqref{eq:phi1} tends to zero exponentially fast
provided that $T_0$ is bounded, or has at most polynomial growth in $z$.
In this case, since by differentiation of \eqref{eq:phi0} the derivative of the interface profile $\partial_z\varphi_0$ belongs to the kernel of the self-adjoint operator 
$\partial_z^2-W''(\varphi_0)$  in $L^2(\R)$ (with domain $H^2(\R)$), 
a necessary condition for \eqref{eq:phi1} to have a solution $\varphi_1(x,t,\cdot)\in H^2(\R)$ is
$$(\overline\alpha\,\partial_t d_0\,-\,\Delta d_0)\,\int_{-\infty}^{+\infty}(\partial_z\varphi_0)^2\,\dif z\,=\,
\overline{\gamma}\,\int_{-\infty}^{+\infty}T_0\,\partial_z \nu(\varphi_0)\,\dif z\,.$$
Defining the `interface temperature' by
$$\langle T_0\rangle\,:=\,\int_{-\infty}^{+\infty}T_0\,\partial_z \nu(\varphi_0)\,\dif z\,=\,\pm\,\frac{\int_{-\infty}^{+\infty}T_0\,\partial_z \nu(\varphi_0)\,\dif z}{\int_{-\infty}^{+\infty}\partial_z \nu(\varphi_0)\,\dif z}\,,$$
and the (nondimensional) surface tension by 
$\sigma_0:=\int_{-\infty}^{+\infty}(\partial_z\varphi_0)^2\,\dif z$,
the previous relation may be seen as a (nondimensional) generalized Gibbs-Thomson condition:
\begin{equation}\label{eq:GB}
\sigma_0\,(\overline\alpha\,\partial_t d_0\,-\,\Delta d_0)\,=\,\pm\,
\overline{\gamma}\,\langle T_0\rangle\,.
\end{equation}
Here above the $\pm$ sign is merely a shorthand for  $[\nu(\varphi_0)]_{-\infty}^{+\infty}$, which equals
$+1$ if $N^0(x,t):=\nabla d_0(x,t)$ points to the liquid phase (or $-1$ of $N^0$ points to the solid phase).
Recalling that $v_0:=-\partial_t d_0$ is the normal speed of $\Gamma^0$ and $H_0:=-\Delta d_0$ is the sum of principal curvatures of $\Gamma^0$, we can indeed identify \eqref{eq:GB} with the usual condition in generalized Stefan models (see again Fig.~1 in \cite{CaginalpChen}).

As regards the second equation in \eqref{eq:systext}, the only term of order zero in $\varepsilon$ is 
$\overline\delta\,\partial_z^2T$. Nevertheless,  we may add to that equation 
a function of the form $$h^\varepsilon(x,t)\,\rho(z)\,(d^\varepsilon(x,t)-\varepsilon z)\,,$$
which obviously vanishes at $z=d^\varepsilon(x,t)/\varepsilon$, with $\rho$ smooth and compactly supported in $z$ 
(so that the term $\varepsilon z$ is at most of the order of $\varepsilon$)
and 
$h^\varepsilon(x,t) \sim \sum_{i=0}^{\infty} \varepsilon^i\,h_i(x,t)$. More precisely, we shall
assume, similarly as in \cite{CaginalpChen}, that $\rho=\partial_z^2\eta$ with $\eta$ such that  $\eta\equiv 0$ on $(-\infty,-1]$, $\eta\equiv 1$ on $[1,+\infty)$, and $\eta'>0$ on $(-1,1)$.  Then the zeroth order equation becomes
$$\overline\delta\,\partial_z^2T_0\,=\,h_0(x,t)\,d_0(x,t)\,\partial_z^2\eta\,,$$
which necessarily yields, if $T_0$ is sought bounded in $z$, that
$$\overline\delta\,T_0(x,t,z)\,-\,h_0(x,t)\,d_0(x,t)\,\eta(z)\,=:\,\overline{\delta}\,{T_0}^-(x,t)\,,$$
a function of $(x,t)$ alone.
Since $\eta(z)=0$ for $z<-1$, this is a consistent notation in that
$T_0^-(x,t)=\lim_{z\to-\infty}T_0(x,t,z)$. Moreover, since $\eta(z)=1$ for $z>1$, we have
$$h_0(x,t)\,d_0(x,t)\,=\,\overline\delta\,(T_0^+(x,t)\,-\,T_0^-(x,t))\,,$$
where $T_0^+(x,t):=\lim_{z\to+\infty}T_0(x,t,z)$.
This shows in particular that for $h_0$ to be smooth, $T_0^+$ and $T_0^-$ must coincide on the zero level set of $d_0$, namely on $\Gamma^0$.  Conversely, if $T_0^+$ and $T_0^-$ are smooth functions coinciding on $\Gamma^0$,
we can define
$$h_0(x,t):=\left\{\begin{array}{ll}
\overline\delta\,\dfrac{T_0^+(x,t)\,-\,T_0^-(x,t)}{d_0(x,t)}\,,& x\notin \Gamma^0(t)\,,\\ [20pt]
\overline\delta\,\left[\nabla T_0\cdot N^0\right](x,t)\,, & x\in \Gamma^0(t)\,,
\end{array}\right.$$ 
where $N^0=\nabla d_0$ (as before) and the `jump' notation $\left[\nabla T_0\cdot N^0\right]$ merely stands for 
$\nabla T_0^+\cdot N^0\,-\,\nabla T_0^-\cdot N^0$.
Then,
$$T_0(x,t,z)\,:=\,\overline{\delta}^{-1}\,h_0(x,t)\,d_0(x,t)\,\eta(z)\,+\,{T_0}^-(x,t)$$
is independent of $z$ for $x\in \Gamma^0(t)$, and more precisely,
$$T_0(x,t,z)\,=\,{T_0}^-(x,t)\,=\,{T_0}^+(x,t)\,=\,\pm\,\langle T_0\rangle (x,t)\,,$$
where again $\pm=[\nu(\varphi_0)]_{-\infty}^{+\infty}$.

Now, the next order terms in the asymptotic expansion of the second equation in  \eqref{eq:systext} supplemented with the term $h^\varepsilon(x,t)\,(d^\varepsilon(x,t)-\varepsilon z)\,\partial_z^2\eta(z)$
will enable us to find a necessary relation between the heat flux $\left[\nabla T_0\cdot N^0\right]$, 
the interface temperature $\langle T_0\rangle +\theta$, and the velocity $v^0=-\partial_td_0$ of $\Gamma^0$.
As a matter of fact, retaining only the terms of order one, we get
$$\begin{array}{l}
\overline\delta\,\partial_z^2T_1\,=\,\begin{array}[t]{l}(h_1\,d_0\,+h_0\,d_1)\,\partial_z^2\eta\,-\,h_0\,z\,\partial_z^2\eta
\\ [10pt] \,+\,
(\overline\beta\partial_td_0-\overline\delta\Delta d_0)\partial_zT_0\,-\,
2\overline\delta\nabla d_0\cdot \nabla\partial_zT_0\\ [10pt]
 +\,\overline\gamma(T_0+\theta)\,(\partial_td_0)\,\partial_z\nu(\varphi_0)\,-\,2\,
\overline\alpha\,(\partial_td_0)^2(\partial_z\varphi_0)^2\,.\end{array}
\end{array}$$
A necessary condition for $T_1$ to be bounded in $z$ is that the integral from $z=-\infty$ to $z=+\infty$ of the right-hand side above equals zero.
The contribution of the first row to the integral is just $h_0$, because
$\int \partial_z^2\eta=0$ and $\int z\partial_z^2\eta=-1$ by definition of $\eta$.
The next term does not contribute if we restrict to $x\in \Gamma^0(t)$ because then $T_0^+(x,t)=T_0^-(x,t)$.
Recalling that $h_0=\overline\delta\,\left[\nabla T_0\cdot N^0\right]$, 
$T_0=\pm \langle T_0\rangle$ on $\Gamma^0$,
and that we have defined
\begin{equation}\label{eq:sigma0}\sigma_0\,=\,
\int_{-\infty}^{+\infty}(\partial_z\varphi_0)^2\,,
\end{equation}
we finally arrive at
\begin{equation}
\label{eq:nablaTbis}
\overline\delta \,\left[\nabla T_0\cdot N^0\right]\,=\,\pm\,
\overline\gamma \,(\langle T_0\rangle + \theta)\,(\partial_td_0)\,-\,2\,\overline\alpha\,\sigma_0\,(\partial_td_0)^2\,.\end{equation}
Observe that the roughly obtained relation \eqref{eq:nablaT} may be seen as an approximation of 
\eqref{eq:nablaTbis} when the velocity $v^0=-\partial_td_0$ of $\Gamma^0$ is small enough.

To summarize, the sharp interface limit of \eqref{eq:systsimpter} is expected to be (rigorous justification will be addressed elsewhere) the generalized Stefan problem
consisting of the heat equation for $T$ outside $\Gamma^0$ together with the conditions \eqref{eq:GB},
\eqref{eq:nablaTbis} on $\Gamma^0$.
If for instance $N=N^0$ points to the liquid phase, this problem reads
\begin{equation}\label{eq:genStefan}
\left\{\begin{array}{ll}\overline\beta\, \partial_t T\,=\,\overline\delta \,\Delta T & \mbox{outside } \Gamma\,,\\
\sigma_0\,(H\,-\,\overline\alpha\,v)\,=\,\overline{\gamma}\,T & \mbox{on } \Gamma\,,\\
\overline\delta \,\left[\nabla T\cdot N\right]\,=\,
-\,\overline\gamma \,(T + \theta)\,v\,-\,2\,\overline\alpha\,\sigma_0\,v^2 & \mbox{on } \Gamma\,,
\end{array}\right.
\end{equation}
with $\sigma_0$ defined in \eqref{eq:sigma0} where $\varphi_0$ is solution 
of \eqref{eq:phi0} and tends to $0$ at $-\infty$ and $1$ at $+\infty$, $v$ the normal velocity of $\Gamma=\Gamma^0$, 
and $H$ the sum of principal curvatures of $\Gamma$.

\subsection{Back to physical variables}\label{ss:backphys}
The sharp interface model \eqref{eq:genStefan}
is non-dimensional. It is of course important from the physical point of view to go back to physical quantities.

Let us start with the first equation in \eqref{eq:genStefan}, which by definition of $\overline\beta$ and $\overline\delta$ also reads
$$\Peclet\,\partial_t T= \nabla T\,.$$
Remembering that $T$, $t$, and $x$ respectively stand for
$$\widetilde{T}=\frac{T-T_e}{\delta T}\,,\quad \widetilde{t}=\frac{t}{t_0}\,,\;\mbox{and }\; \widetilde{x}=\frac{x}{L}\,,$$
(where the tilda are those of Section \ref{ss:nondim} and not those of Section \ref{ss:deras})
by definition of the Peclet number $\Peclet={\rho C_p L^2}/({\cond t_0})$ we recover the expected heat equation
$$\rho\,C_p \partial_t T\,=\,\cond \,\Delta T\,.$$

As to the last equation in  \eqref{eq:genStefan}, it actually reads
$$
\frac{1}{\Peclet}\,\left[\nabla_{\widetilde{x}} \widetilde{T}\cdot N\right]\,=\,
-\,\frac{1}{\Stefan\theta} \,(\widetilde{T} + \theta)\,\widetilde{v}\,-\,2\,\frac{\beta\,\sigma_0}{\alpha\,\varepsilon}\,\widetilde{v}^2$$
with $\widetilde{v}=\partial_{\widetilde{t}}d/L=(t_0/L)v$ if $v$ denotes the physical velocity of the interface. 
Before going further, let us comment on 
$\sigma_0=\int_{-\infty}^{+\infty}(\partial_z\varphi_0)^2\,\dif z$, where $\varphi_0$ is by definition (see Eq.~\eqref{eq:phi0}, having in mind that  $W$ stands for $\widetilde{W}=(h/\sigma)\,W$) solution of the differential equation
$$\partial_z^2\varphi_0\,=\,\widetilde{W}'(\varphi_0)\,.$$
An obvious integrating factor is $\partial_z\varphi_0$, and since $W(\varphi_0)$ vanishes at $\pm \infty$, we have
$(\partial_z\varphi_0)^2=2(h/\sigma)\,W(\varphi_0)$. This implies that 
$$\sigma_0=\frac{2}{\sigma}\int_{-\infty}^{+\infty}W(\varphi_0)\,h\,\dif z\,.$$
Recalling the meaning of the $z$ variable, which scales as the actual distance to the interface over $\varepsilon L=h$,
we can identify the integral $2\int_{-\infty}^{+\infty}W(\varphi_0)\,h\,\dif z$ with the physical \emph{surface tension} $\sigma$, and therefore set $\sigma_0=1$.
Substituting all the other non-dimensional parameters $\Peclet$, $\Stefan$, $\theta$, $\beta$, $\alpha$ and $\varepsilon$ by their expressions in terms of physical quantities, we get in turn the physical jump condition
$$\frac{k}{\rho}\, \left[\nabla {T}\cdot N\right]\,=\,-\,{\mathcal L}\,v\,-\,2\,\frac{v^2}{\kappa\,h}\,.$$
(Recall from \eqref{eq:simplified} that ${\mathcal L}={\mathcal L}_e\,{T}/{T_e}$.)
This is to be compared with the  usual jump condition in Stefan models:
$$\frac{k}{\rho}\, \left[\nabla {T}\cdot N\right]\,=\,-\,{\mathcal L}\,v\,.$$
In particular, we observe that the quadratic correction in the velocity $v$
is negligible if
the velocity $v$ is small compared to ${\mathcal L}\,\kappa\, h$ (which is indeed homogeneous to a velocity).

We finish with the derivation of the generalized Gibbs-Thomson relation.
Recalling that we have set $\sigma_0=1$, that $v$ stands for 
$\widetilde{v}=(t_0/L)v$ and noting that $H$ stands for 
$\widetilde{H}=\Delta_{\widetilde x} \widetilde{d}=L\Delta {d}$, we can rewrite the second equation in  \eqref{eq:genStefan} as 
$$\overline\alpha\,\frac{t_0}{L}\,v\,-\,L\,H\,=\,-\,\overline{\gamma}\,\widetilde{T}\,,$$
where $H:=\Delta d$ is the actual sum of principal curvatures (homogeneous to the inverse of a distance).
Substituting $\overline\alpha$ and $\overline{\gamma}$ for their expressions, this eventually gives
$$-\,\frac{\rho}{\kappa\,h}\,v\,+\,\sigma\,H\,=\,\rho\,\frac{{\mathcal L}_e}{T_e}\,(T-T_e)\,.$$

Therefore, in physical variables the (generalized Stefan) sharp interface model \eqref{eq:genStefan} reads
\begin{equation}\label{eq:genStefanphys}
\left\{\begin{array}{ll}  \rho\,C_p\, \partial_t T\,=\,\cond \,\Delta T & \mbox{outside } \Gamma\,,\\ [10pt]
\rho\,\dfrac{{\mathcal L}_e}{T_e}\,(T-T_e)\,=\,-\,\dfrac{\rho}{\kappa\,h}\,v\,+\,\sigma\,H & \mbox{on } \Gamma\,,\\ [10pt]
\dfrac{k}{\rho}\, \left[\nabla {T}\cdot N\right]\,=\,-\,{\mathcal L}\,v\,-\,2\,\dfrac{v^2}{\kappa\,h} & \mbox{on } \Gamma\,.
\end{array}\right.
\end{equation}

\section{Well-posedness}\label{s:well}

We now turn to the mathematical analysis of the (non-standard) PDEs system~\eqref{eq:systsimpbis}, which we equivalently rewrite as
\begin{equation}\label{eq:systsimp1}
\left\{
\begin{aligned}
& \hat\alpha \, \partial_t\varphi -\varepsilon^2\Delta \varphi+W'(\varphi)=\gamma \, \nu'(\varphi) \, T, \\
&\hat\beta \, \partial_t T +\gamma \, \theta \, \nu'(\varphi) \, \partial_t\varphi- \delta \, \Delta T=F(\varphi,\Delta\varphi,T),
\end{aligned}
\right.
\end{equation}
with 
\begin{equation}\label{eq:systsimp1F}
F(\varphi,\Delta\varphi,T):=\frac{1}{\hat\alpha}(\ep^2\Delta\varphi-W'(\varphi))^2+\frac{\gamma}{\hat\alpha}(\ep^2\Delta\varphi-W'(\varphi))\nu'(\varphi)T.
\end{equation}
In this system, $\hat\alpha$, $\hat\beta$, $\gamma$, $\delta$, $\varepsilon$, and $\theta$ are fixed, positive parameters,
and the functions $W$, $\nu$ are supposed to be nonnegative and to belong to ${\mathscr C}^3_b(\R)$ (the space of ${\mathscr C}^3$ functions that are bounded as well as their derivatives up to order $3$).
We are going to show that the Initial Boundary Value Problem for \eqref{eq:systsimp1} with suitable initial and boundary data is locally well-posed both in two and three space dimensions.

\subsection{Functional framework and main results}

In what follows,  $\Omega$ is an open, bounded, and regular subset of $\R^d$, $d\in\{2,3\}$.
As boundary conditions 
on $\partial\Omega$ we consider a homogeneous Neumann condition for the order parameter $\varphi$, and a mixed constant Neumann-Dirichlet boundary condition for the temperature:
\begin{equation}\label{boundary:systsimp1}
\frac{\partial \fy}{\partial {\mathrm {n}}} = 0\text{ on $\partial\Omega$}, \qquad \frac{\partial T}{\partial {\mathrm {n}}} = q_{\mathrm{b}} \text{ on $\Gamma$},\qquad T=T_{\mathrm{b}} \text{ on $\partial\Omega\setminus\Gamma$},
\end{equation}
where $\Gamma$ is a given, relatively open subset of $\partial\Omega$, and ${\bf n}$ denotes the unit outward normal to $\partial\Omega$.
We suppose that $q_{\mathrm{b}}$ and $T_{\mathrm{b}}$ are constants, corresponding respectively to the heat flux and to the temperature imposed on the boundary of the domain.
Given $q_{\mathrm{b}}$, $T_{\mathrm{b}}$, we know from \cite[Notes of chapter 8]{GilbargTrudinger} that there exists $\widetilde T\in H^1(\Omega)$ solution of 
\begin{equation}\label{mixedT}
\left\{
\begin{aligned}
\Delta \widetilde T & = 0 \quad \text{in }\Omega,\\
\displaystyle \frac{\partial \widetilde T}{\partial n} & = q_{\mathrm{b}} \quad \text{on }\Gamma,\\
\widetilde T & = T_{\mathrm{b}} \quad \text{on }\partial\Omega\setminus\Gamma,
\end{aligned}
\right.
\end{equation}
in the sense that $\widetilde T-T_{\mathrm{b}} \in H^1_0(\Omega\cup\Gamma)$, 
the closure of ${\mathscr C}^1_0(\Omega\cup\Gamma)$ in $H^1(\Omega)$, 
and 
\begin{equation*}
\int_\Omega \nabla\widetilde T\cdot\nabla \tau dx=\int_\Gamma q_{\mathrm{b}} \tau d{\mathcal{H}^{d-1}}(x)
\end{equation*}
for all $\tau\in H^1_0(\Omega\cup\Gamma)$.
Two solutions $\widetilde T_1$ and $\widetilde T_2$ of the first two equations in~\eqref{mixedT} differ by a constant, so that $\widetilde T$ is unique if $|\partial\Omega\setminus\Gamma|>0$, and unique up to constant in the case of a pure Neumann condition.

For $s\geq 0$, we denote by 
$H^s_{\mathrm {n}}(\Omega)$ the closure of 
\begin{equation*}
\Big\{ \varphi\in \mathcal D(\overline \Omega)~;\quad \frac{\partial \varphi}{\partial n}\Big|_{\partial\Omega} = 0 \;  \Big\}
\end{equation*}
in $H^s(\Omega)$. 
In particular, $H^2_{\mathrm {n}}(\Omega)$ is merely the set of functions $\varphi\in H^2(\Omega)$ that satisfy the homogeneous Neumann boundary condition $\partial_{\mathrm {n}}\varphi=0$ on $\partial\Omega$. 
When no confusion can occur we shall just write $L^2$ for $L^2(\Omega)$, $H^s_{\mathrm {n}}$ for $H^s_{\mathrm {n}}(\Omega)$, and 
$H^1_0$ for $H^1_0(\Omega\cup \Gamma)$.

\begin{definition}[Weak solution]
\label{definition-weaksolution}
For $\fy_0\in H^2_{\mathrm {n}}$, $T_0\in H^1_0$, and $t^*\in (0,+\infty]$,
we say that $(\fy, T)$ is a \emph{weak solution} of \eqref{eq:systsimp1} on $[0,t^*)$ with initial data $\fy_0$, $T_0$, and boundary conditions~\eqref{boundary:systsimp1}, if
\begin{equation*}
\begin{aligned}
& \fy \in {\mathscr C}_b ([0,t^*);H^2_{\mathrm {n}}) \cap L^2_{\text{loc}}(0,t^*;H^3_{\mathrm {n}}), 
\\
& T-\widetilde T \in   {\mathscr C}_b([0,t^*);L^2) \cap L^2_{\text{loc}}(0,t^*;H^1_0),
\end{aligned}
\end{equation*}
with $\fy$ and $T$ satisfying $\fy \big|_{t=0} = \fy_0, \quad T \big|_{t=0} = T_0$, and
\begin{itemize}
\item the first equation in \eqref{eq:systsimp1} in the sense that
\begin{equation}
 \hat\alpha\,\int_{0}^{t^*}\!\!\!\!\int_{\Omega} \fy \partial_t\zeta  = \int_{0}^{t^*}\!\!\!\!\int_{\Omega} W'(\varphi)\,\zeta - \gamma \int_{0}^{t^*}\!\!\!\!\int_{\Omega} \nu'(\varphi)\,T\, \zeta - \varepsilon^2 \int_{0}^{t^*}\!\!\!\!\int_{\Omega} \Delta \varphi\, \zeta
\label{eq:weaksolvarphi}\end{equation}
for all $\zeta \in {\mathscr C}^1_0((0,t^*);L^2)$,
 \item  the second equation in \eqref{eq:systsimp1}  in the sense that 
\begin{multline}
\int_{0}^{t^*}\!\!\!\!\int_{\Omega}  (\hat\beta T+ \gamma\theta\nu(\varphi) )\,\partial_t\tau  =  \,\delta \int_{0}^{t^*}\!\!\!\!\int_{\Omega} \nabla T \cdot \nabla \tau 
\,
-\,\int_{0}^{t^*}\!\!\!\!\int_\Omega F(\varphi,\Delta\varphi, T)\,\tau
\label{eq:weaksolT}\end{multline}
for all $\tau \in {\mathscr C}^1_0((0,t^*);L^2)\cap {\mathscr C}((0,t^*);H^1_0)$,
where $F$ is defined as in \eqref{eq:systsimp1F}.
\end{itemize}
\end{definition}

Note that, according to this definition and by the Sobolev embeddings
$H^2\hookrightarrow {\mathscr C}^0$ and $H^1 \hookrightarrow L^6$ (both valid in space dimension $d\leq 3$), a weak solution $(\fy, T)$ is such that $\varphi\in {\mathscr C}_b ([0,t^*)\times \Omega)$, and $\Delta\varphi\in L^6(\Omega)$ at almost all times in $[0,t^*)$.
This gives  sense in particular to the last integral in~\eqref{eq:weaksolT} if we also
note that $\tau \in L^6(\Omega)$ at all times in $(0,t^*)$ if $\tau \in{\mathscr C}((0,t^*);H^1_0)$: indeed, examining all the terms in the product 
$F(\varphi,\Delta\varphi, T)\,\tau$ we see, using that $W'$ and $\nu'$ are bounded, that we `only' need
$\tau$, $T\,\tau$, $\Delta\varphi\,T\,\tau$, $\Delta \varphi\,\tau$, and $(\Delta\varphi)^2\,\tau$ being integrable, which is certainly the case on a bounded domain when $T\in L^2$, $\tau\in L^6$, and $\Delta\varphi\in L^6$.

\begin{theorem}[Existence of weak solutions]
\label{theorem-localweaksolution}
For all $\fy_0\in H^2_{\mathrm {n}}$, $T_0\in H^1_0$,
there exist $t^*>0$ and a weak solution $(\fy, T)$  of \eqref{eq:systsimp1} on $[0,t^*)$ with initial data $\fy_0$, $T_0$, and boundary conditions~\eqref{boundary:systsimp1}, in the sense of Definition~\ref{definition-weaksolution}.
\end{theorem}
\begin{theorem}[Continuous dependence on the data]
\label{theorem-unicitesolution}
\begin{itemize}
\item[1.] Given  $\fy_0\in H^2_{\mathrm {n}}$, $T_0\in H^1_0$, there exists at most one weak solution $(\fy, T)$ to~\eqref{eq:systsimp1}--\eqref{boundary:systsimp1} with initial datum $(\fy_0,T_0)$.
\item[2.] If $(\varphi_i,T_i)$, $i\in \{1,2\}$, are two weak solutions to~\eqref{eq:systsimp1}--\eqref{boundary:systsimp1} both defined on $[0,t^*]$, and with initial data $(\fy_{i,0},T_{i,0})$, $i=1,2$, respectively, then there exists a constant $C>0$ depending only on the norms 
\begin{equation*}
\begin{aligned}
& \|T_i\|_{L^\infty(0,t^*;H^1(\Omega))},\quad
\|T_i\|_{L^2(0,t^*;L^2(\Omega))},\\
& \|\varphi_i\|_{L^\infty(0,t^*;H^2(\Omega))},\quad
\|\varphi_i\|_{L^2(0,t^*;H^3(\Omega))},
\end{aligned}
\qquad i\in \{1,2\}
\end{equation*}
such that, for $0\leq t\leq t^*$, 
\begin{multline}
\|T_1(t)-T_2(t)\|_{L^2(\Omega)}^2+\|\varphi_{1}(t)-\varphi_{2}(t)\|_{H^2(\Omega)}^2\\
 \leq C\left(\|T_{1,0}-T_{2,0}\|_{L^2(\Omega)}^2+\|\varphi_{1,0}-\varphi_{2,0}\|_{H^2(\Omega)}^2\right).
\label{eq:continuousdependencedata}\end{multline}
In particular, the weak solution to~\eqref{eq:systsimp1}--\eqref{boundary:systsimp1} depends continuously on the data.
\end{itemize}
\end{theorem}

Unsurprisingly, the proof of Theorems \ref{theorem-localweaksolution} and~ \ref{theorem-unicitesolution}
relies on \emph{a priori} estimates. It is to be noted though that we shall use other quantities than the total energy
$$E(T,\varphi)\,:=\,\int_\Omega\, \big(\hat\beta\, T\,+\,\gamma\theta\,\nu(\varphi)\,+\,W(\varphi)\,+\,\tfrac{1}{2}\,\varepsilon^2\,|\nabla\varphi|^2\big)\,.$$
Indeed, if we do have the conservation of $E$, thanks to \eqref{eq:evolutionenondim} (where $e=\beta E$), at least when 
$T$ satisfies a homogeneous Neumann condition on the whole boundary $\partial\Omega$
(that is, for $q_{\mathrm{b}}=0$ and $\Gamma=\partial\Omega$), this is obviously not enough to control the $L^2$ norm of $T$.
Rather, we shall use
\begin{equation}\label{eq:E0}
E_0(T,\varphi)\,:=\,\int_\Omega\, \big(\frac{\hat\beta}{2\theta}\, T^2\,+\,W(\varphi)\,+\,\tfrac{1}{2}\,\varepsilon^2\,|\nabla\varphi|^2\big)\,.
\end{equation}
Compared to $E$, the interest of $E_0$ is that it is quadratic in $T$, and it satisfies the identity 
$${\theta}\,\frac{\dif }{\dif t} E_0(T,\varphi)\,+\,{\delta}\,\|\nabla T\|^2_{L^2}\,+\,\hat\alpha\,{\theta}\,\|\partial_t\varphi\|^2_{L^2}\,=\,\int_\Omega 
F(\varphi,\Delta\varphi,T)\,T\,$$
along solutions of \eqref{eq:systsimp1} (and \eqref{boundary:systsimp1} with $q_{\mathrm{b}}=0$, $T_{\mathrm{b}}=0$).
Of course, because of the right-hand side, this is not fully satisfactory  and we shall need another quantity to control the $L^2$ norm of $\Delta\varphi$.

\begin{remark}
If we replace $F$ by zero in \eqref{eq:systsimp1}, we are left with the Caginalp-like model
\begin{equation}\label{eq:Caginalp}
\left\{
\begin{aligned}
& \hat\alpha \, \partial_t\varphi -\varepsilon^2\Delta \varphi+W'(\varphi)=\gamma \, \nu'(\varphi) \, T, \\
&\hat\beta \, \partial_t T +\gamma \, \theta \, \nu'(\varphi) \, \partial_t\varphi- \delta \, \Delta T=0,
\end{aligned}
\right.
\end{equation}
for which we have the much nicer identity
$$\theta\,\frac{\dif }{\dif t} E_0(T,\varphi)\,+\,{\delta}\,\|\nabla T\|^2_{L^2}\,+\,\hat\alpha\,{\theta}\,\|\partial_t\varphi\|^2_{L^2}\,=\,0\,.$$
\end{remark}

\begin{remark} By an adaptation of our \emph{a priori} estimates, we can show in addition that, among the stationary solutions to~\eqref{eq:systsimp1}, those corresponding to single-phase states, {\it i.e.} with $\varphi\equiv 0$ or $\varphi\equiv 1$ and $T=\widetilde T$, are stable. 
\end{remark}

\begin{remark} Since we are interested in asymptotic models, we have chosen on purpose to keep track of the nondimensionalized numbers $\hat\alpha$, $\hat\beta$, $\gamma$, $\delta$, $\varepsilon$, and $\theta$  in our  \emph{a priori} estimates. We shall also pay attention to the occurrence of  $\|\widetilde T\|_{H^1}$, and of the bounds for $W$, $\nu$, and their derivatives.
\end{remark}

\paragraph{Further notations.} All constants depending only on the dimension $d$ and on $\Omega$  will be considered 
harmless, and we shall denote by
\begin{equation*}
A_1\lesssim A_2
\end{equation*}
any inequality $A_1\leq C A_2$ where $C$ is a constant (depending only on $d$ and $\Omega$). 
As already mentioned, $W$ and $\nu$ are supposed to belong to ${\mathscr C}^3_b(\R)$. 
For simplicity, we introduce the notations
\begin{equation*}
\nu'_\infty:=\sup_\R|\nu'|,\quad W'_\infty:=\sup_\R|W'|,
\end{equation*}
and similarly for their second and third order derivatives.
\medskip

\subsection{Existence of solutions}


In this section, we prove Theorem~\ref{theorem-localweaksolution} by means of  a Galerkin approximation.

Let $\{\varphi_i\}_{i\in\N^*}$ be a set of eigen-functions of the Laplacian operator $-\Delta$ with Neumann boundary condition,  $\{\varphi_i\}_{i\in\N^*}$ being a complete orthonormal system in ${H^1}$. Let $\{\overline T_i\}_{i\in\N^*}$ be a set of eigen-functions of $-\Delta$ in $\Omega$ with the boundary conditions
\begin{equation*}
\frac{\partial\overline T_i}{\partial n}=0\mbox{ on }\Gamma,\quad\overline T_i=0\mbox{ on }\partial\Omega\setminus\Gamma,
\end{equation*}
with $\{\overline T_i\}_{i\in\N^*}$ being a complete orthonormal set in ${H^1_0}$. We seek approximate solutions
of \eqref{eq:systsimp1}~-~\eqref{boundary:systsimp1} of the form
\begin{equation*}
\varphi^n(t)=\sum_{i=1}^n a_i(t)\varphi_i,\quad T^n(t):=\overline T^n(t)+\widetilde T\,,\quad  \overline T^n(t)=\sum_{i=1}^n b_i(t)\overline T_i, 
\end{equation*}
where $a_i$ and $b_i$ are  ${\mathscr C}^1$, real-valued functions. Defining
\begin{equation*}
\mathcal V_n=\mathrm{Span}\{\varphi_1,\cdots, \varphi_n\},\quad \mathcal Z_n=\mathrm{Span}\{\overline T_1,\cdots, \overline T_n\},
\end{equation*}
we require that for all $\zeta \in \mathcal V_n$,
\begin{equation}\label{galerkin-fy}
\int_{\Omega} \hat\alpha\,\partial_t\fy^n \, \zeta = -\int_{\Omega} W'(\fy^n)\zeta +\gamma \int_{\Omega} \nu'(\varphi^n) T^n \zeta +\varepsilon^2 \int_{\Omega} \Delta \varphi^n  \zeta,
\end{equation}
and for all $\tau \in \mathcal Z_n$,
\begin{equation}\label{galerkin-T}
\int_{\Omega} \hat\beta\,\partial_t\overline T^n\, \tau = -\delta\int_{\Omega} \nabla \overline T^n \cdot \nabla \tau 
-\gamma\theta\int_\Omega\nu'(\varphi^n)\partial_t\varphi^n \tau
+\int_\Omega F^n \tau,
\end{equation}
where
\begin{equation*}
F^n:=\frac{1}{\hat\alpha}(\ep^2\Delta\varphi^n-W'(\varphi^n))^2+\frac{\gamma}{\hat\alpha}(\ep^2\Delta\varphi^n-W'(\varphi^n))\nu'(\varphi^n)T^n\,,
\end{equation*}
together with the initial conditions 
$$\fy^n(0)=\mathcal P_{\mathcal V_n}(\fy_0)\,,\quad \overline T^n(0)=\mathcal P_{\mathcal Z_n}(T_0-\widetilde T)\,,$$
where for any subspace $Y$, $\mathcal P_Y$ denotes the orthogonal projection onto $Y$.
\medskip

Denoting by $\mathbf a$ and $\mathbf b$  the vector-valued functions 
of components $a_i$ and $b_i$ respectively, taking $\zeta=\varphi_i$ in~\eqref{galerkin-fy} and $\tau=\overline T_i$ in~\eqref{galerkin-T} for $i=1,\ldots,n$, we obtain ordinary differential equations of the form
\begin{equation}
\frac{\dif {\mathbf a}}{\dif t}=\Phi(\mathbf a,\mathbf b)\,,\quad \frac{\dif {\mathbf b}}{\dif t}=\mathcal{T}(\mathbf a,\mathbf b,\Phi(\mathbf a,\mathbf b))\,,
\label{eq:ODE}
\end{equation}
with $\Phi$ and $\mathcal{T}$ at least ${\mathscr C}^2$ on $\R^{2n}$ and $\R^{3n}$ respectively (since 
$W$ and $\nu$ are ${\mathscr C}^3$). Therefore, the Cauchy-Lipschitz theorem ensures the existence and uniqueness on some maximal time interval $[0,t_n)$, $t_n\in (0,+\infty]$, of a solution 
$(\mathbf a,\mathbf b)$ with prescribed initial data. 
\medskip

The next step is to derive some estimates in order to show that the times sequence~$(t_n)$ is bounded from below by a positive time and that the sequences~$(\fy^n)$, $(\overline T^n)$ are bounded in the appropriate functional spaces.
To simplify the notations we will drop the 
superscript~$n$. 

\medskip

\paragraph{Energy estimate.} 
Recalling the definition of $E_0$ in \eqref{eq:E0}, and taking a combination of~\eqref{galerkin-fy} with $\zeta=\varphi_n$ and~\eqref{galerkin-T} with $\tau=\overline T_n$, we get the identity
\begin{equation}
\theta\,\displaystyle \frac{\dif\;}{\dif t}E_0(\overline T,\varphi)+{\delta}\,\|\nabla \overline T\|^2_{L^2}\,+\,\hat\alpha\,{\theta}\,\|\partial_t\varphi\|^2_{L^2}=\int_\Omega F\overline T+\gamma\theta\,\int_\Omega \nu'(\varphi)\partial_t\varphi \,\widetilde T.
\label{eq:energyCaginalpcoeff}\end{equation}
By the elementary inequality 
\begin{equation}\label{eq:elem}
ab\leq \lambda a^2+\frac{b^2}{4\lambda}\,,\quad \forall a, b \in \R\,,\;\forall \lambda>0\,,
\end{equation}
we have
\begin{equation*}
\gamma\int_\Omega \nu'(\varphi)\partial_t\varphi\, \widetilde T\leq \frac{\hat\alpha}{4}\,\|\partial_t\varphi\|^2_{L^2}\,+
\frac{\gamma^2}{\hat\alpha}\,|\nu'_\infty|^2\,\|\widetilde T\\|_{L^2}^2.
\end{equation*}
Recalling that the notation $\hat\alpha$ actually means $1/\alpha$, and introducing the new simplifying notation
\begin{equation}\label{eq:mu}
\mu:=\gamma^2\,|\nu'_\infty|^2\,,
\end{equation}
we thus obtain from~\eqref{eq:energyCaginalpcoeff} that
\begin{equation}
\displaystyle\frac{\dif}{\dif t}E_0(\overline T,\varphi)+
\delta\hat\theta\,\|\nabla \overline T\|^2_{L^2}+\tfrac{3}{4}\hat\alpha\,\|\partial_t\varphi\|^2_{L^2}\leq
\hat\theta\int_\Omega F\overline  T+
\alpha\,\mu\,\|\widetilde T\|_{L^2}^2,
\label{eq:energyCaginalpcoeff2}\end{equation}
where we have used the notation $\hat\theta$ for $1/\theta$. We shall use the same convention repeatedly for other positive parameters in what follows.

The integral of  $F\overline T$ in~\eqref{eq:energyCaginalpcoeff2} 
involves trilinear terms in $\Delta \varphi$ and $T$. We can get an estimate on 
them by using the following, higher order estimate.
\medskip

\begin{lemma}[second order estimate] Let us define the modified energy
\begin{equation*}
E_1(
T,\varphi)=E_0(
T,\varphi)+\tfrac12
{\ep^2\delta\hat\alpha}{\hat\mu\hat\theta}\,\int_\Omega |\Delta\varphi|^2\,,
\end{equation*}
and (to enforce estimates of powers of $E_1$)
\begin{equation*}
E_1^*(
T,\varphi)=\max(1,E_1(
T,\varphi)).
\end{equation*}
Let the constants $A^0$, $B^0$, $C^0$, $D^0$ be given by
\begin{align*}
A^0:=\,&
\delta
\omega\hat\mu \hat\theta \hat\varepsilon^2,\quad
D^0:=\,
\max 
(1,{\theta}{\beta}
)\,{\theta}{\beta}\,(1+ 
\mu\alpha\theta\hat\delta
)\,\iota\mu\hat\ep^{4} (1+ 
\iota\mu\hat\ep^4
),\\
B^0:=\,&
\hat\beta\hat\theta\|\widetilde T\|_{L^2}^2,\quad
C^0:=\,
(\mu\alpha
+
\delta\hat\theta)
 \|\widetilde T\|^2_{H^1},
\end{align*}
where
$$\omega:=|W''_\infty|^2\,,\quad\mu:=\gamma^2\,|\nu'_\infty|^2\,,\quad \iota:=\left|\dfrac{\nu''_\infty}{\nu'_\infty}\right|^2\,.$$
Then the first $n$ modes $(T=\overline T+\widetilde{T},\varphi)=(T^n=\overline T^n+\widetilde{T},\varphi^n)$ solutions to~\eqref{galerkin-fy}-\eqref{galerkin-T} satisfy the refined energy estimate
\begin{multline}
\displaystyle \frac{\dif\;}{\dif t}E_1(\overline T,\varphi)+
\delta\hat\theta\|\nabla\overline  T\|_{L^2}^2+\hat\alpha\|\partial_t\varphi\|_{L^2}^2+
\ep^4\delta\hat\mu\hat\theta\,\|\nabla\Delta\varphi\|_{L^2}^2\\
-
\hat\theta\int_\Omega F(\varphi,\Delta\varphi,T)\,\overline  T \lesssim A^0 E_1(\overline T,\varphi)+ D^0(B^0+E_1^*(\overline T,\varphi))^3+C^0,
\label{eq:totenergyestimate0}\end{multline}
for all $t\in(0,t_n)$.
\label{lem:E1}\end{lemma}

{\bf Proof:} Let us first write
\begin{equation*}
E_1(\overline T,\varphi)=E_0(\overline T,\varphi)+z\hat\alpha\int_\Omega |\Delta\varphi|^2,
\end{equation*}
with the constant~$z$ to be determined in the course of the proof. 
We apply~\eqref{galerkin-fy} with~$\Delta^2 \fy$ as a test function.
Since we use a Galerkin approximation built on the eigen-functions of the operator $-\Delta$ with Neumann boundary condition, there is no boundary term in the next integrations by parts.
We obtain
\begin{equation}\label{estimate11}
\tfrac{1}{2}
 \hat\alpha \frac{\dif}{\dif t} \|\Delta \fy\|_{L^2}^2
 +\ep^2 \left\| \nabla \Delta \fy \right\|_{L^2}^2
= - \underbrace{\int_{\Omega} W'(\fy)\Delta^2 \fy}_{J_1}  + \underbrace{\gamma \int_{\Omega} \nu'(\fy) T\, \Delta^2 \fy}_{J_2}.
\end{equation}
Using an integration by parts, we get an estimate on~$J_1$ because
\begin{equation*}
 \left| \int_{\Omega} W''(\fy) \nabla \fy \cdot \nabla \Delta \fy \right|\,-\,
\dfrac{\ep^2}{4} \left\| \nabla \Delta \fy \right\|_{L^2}^2 
\lesssim 
 \frac{|W''_\infty|^2}{\varepsilon^2} \|\nabla \fy\|_{L^2}^2
\end{equation*}
thanks to the  
inequality in \eqref{eq:elem}.
Another integration by parts gives
\begin{equation*}
J_2
= - \gamma  \int_{\Omega} \nu'(\fy) \nabla T \cdot \nabla \Delta \fy 
  - \gamma  \int_{\Omega} \nu''(\fy) T\, \nabla \fy \cdot \nabla \Delta \fy.
\end{equation*}
By the Cauchy-Schwarz inequality, we thus have
\begin{equation*}
|J_2|
\leq \gamma \nu'_\infty \|\nabla T\|_{L^2} \|\nabla \Delta \fy\|_{L^2}
+ \gamma\nu''_\infty \|T\|_{L^2} \|\nabla \fy\|_{L^\infty} \|\nabla \Delta \fy\|_{L^2},
\end{equation*}
hence by using \eqref{eq:elem} again,
\begin{equation*}
|J_2|\,-\, 
\dfrac{\ep^2}{8}\,\|\nabla \Delta \fy\|_{L^2}^2
\lesssim 
\left|\dfrac{\gamma\nu'_\infty}{\ep}\right|^2  \|\nabla T\|^2_{L^2}
+ \gamma \nu''_\infty \|T\|_{L^2} \|\nabla \fy\|_{L^\infty} \|\nabla \Delta \fy\|_{L^2}.
\end{equation*}
To control the $L^\infty$ norm of $\nabla\varphi$, we can apply to $u:=\nabla\varphi$ Agmon's inequality \begin{equation}
\|u \|_{L^\infty}^2 \lesssim \|u\|_{H^1}\, \|u \|_{H^2}\,,
\label{eq:Agmon}\end{equation}
which is valid in dimension $d\leq 3$ -~in the two dimensional case, it just  follows from the Sobolev embedding $H^{3/2}
\hookrightarrow L^\infty
$ and the interpolation between $H^1$ and $H^2$, showing that
$H^{3/2}=[H^1,H^2]_{1/2}$. Using that $\|u\|_{H^2}\leq \|u\|_{L^2}+\|\Delta u\|_{L^2}$, this gives
\begin{equation*}
\|\nabla \fy \|_{L^\infty} \lesssim \|\nabla \fy\|_{H^1}^{1/2} ( \| \nabla \fy \|_{L^2}^{1/2} + \| \nabla \Delta \fy \|_{L^2}^{1/2} ),
\end{equation*}
from which we deduce that 
\begin{equation*}
\|T\|_{L^2} \|\nabla \fy\|_{L^\infty} \|\nabla \Delta \fy\|_{L^2} \lesssim \|T\|_{L^2}( \|\nabla \fy\|_{H^1} \| \nabla \Delta \fy \|_{L^2} + \|\nabla \fy\|_{H^1}^{1/2} \| \nabla \Delta \fy \|_{L^2}^{3/2}).
\end{equation*}
Then, using once more \eqref{eq:elem} together with its more general version, Young's inequality
\begin{equation}\label{eq:Young}
ab\leq \lambda \frac{a^p}{p}+\lambda^{-q/p}\,\frac{b^q}{q},\quad  \forall a, b \in \R\,,\;\forall \lambda>0\,,\quad \, \forall p,q>0\,,\;\frac{1}{p}+\frac{1}{q}=1,
\end{equation} 
with $p=4/3$,
 we obtain
\begin{equation*}
|J_2|\,-\, 
\dfrac{\ep^2}{4}\,\|\nabla \Delta \fy\|_{L^2}^2
\lesssim 
\mu\hat\ep^2\|\nabla T\|^2_{L^2}
+
\iota\mu\hat\ep^2\|T\|^2_{L^2}\|\nabla\varphi\|_{H^1}^2
+
\iota^2\mu^2\hat\ep^6\|T\|^4_{L^2}\|\nabla\varphi\|_{H^1}^2,
\end{equation*}
where we have used the announced shorthands
$$\mu=\gamma^2\,|\nu'_\infty|^2\,,\quad \iota=\left|\dfrac{\nu''_\infty}{\nu'_\infty}\right|^2\,.$$
Now, since $W\geq 0$ by assumption, we have
$$\|\nabla\fy\|_{L^2}^2\,\leq\,
2\hat\ep^2\,E_0(T,\varphi)
\,,\quad
\|T\|_{L^2}^2\,\leq\,
2\theta\beta\,E_0(T,\varphi)\,
\,,$$
and $E_0\leq E_1^*$,
so that we can rewrite the estimate of $J_2$ as
\begin{equation*}
|J_2| -\, 
\dfrac{\ep^2}{4}\,\|\nabla \Delta \fy\|_{L^2}^2
\lesssim 
\mu\hat\ep^2 \|\nabla T\|^2_{L^2}
+
\zeta\,E_1^*(T,\varphi)^3\,,
\end{equation*}
with 
$$\zeta:=\max \left(1,{2\theta}{\beta}\right)\,{2\theta}{\beta}\,\left({2}{\hat\ep^2}+{\hat{z}\alpha}\right)\,
\iota\mu\hat\ep^2(1+\iota\mu\hat\ep^4)
\,.$$
Recalling the estimate of $J_1$ obtained at the beginning, in which we can bound $\|\nabla\varphi\|_{L^2}^2$ by $2\hat\ep^2\,E_1(\overline{T},\varphi)$,
we find that 
\begin{multline*}
\hat\alpha\, \frac{\dif}{\dif t}  \|\Delta \fy\|_{L^2}^2
 +\,{\ep^2} \left\| \nabla \Delta \fy \right\|_{L^2}^2
\lesssim 
\omega\hat\ep^4
E_1(\overline T,\varphi)
+ 
\mu\hat\ep^2 \|\nabla T\|^2_{L^2}
+
\zeta\,E_1^*(T,\varphi)^3,
\end{multline*}
where we have used the other shorthand
$$\omega=|W''_\infty|^2\,.$$
Finally, noting that $\|T\|_{L^2}\leq \|\overline T\|_{L^2}+\|\widetilde T\|_{L^2}$, and similarly for the derivatives of~$T$, we obtain
\begin{multline}\label{estimate_fy1}
 \hat\alpha \,\frac{\dif}{\dif t} \|\Delta \fy\|_{L^2}^2
 +{\ep^2} \left\| \nabla \Delta \fy \right\|_{L^2}^2
\lesssim  
\mu\hat\ep^2 \|\nabla \overline T\|^2_{L^2}
+ 
\mu\hat\ep^2 \|\nabla \widetilde T\|^2_{L^2}\\
+
\omega\hat\ep^4E_1(\overline T,\varphi)+
\zeta\,\left(
\hat\beta\hat\theta\|\widetilde T\|_{L^2}^2+E_1^*(\overline T,\varphi)\right)^3.
\end{multline}
Multiplying this inequality by $z$ and adding it with \eqref{eq:energyCaginalpcoeff2}, we can absorb  $\|\nabla \overline T\|^2_{L^2}$ on the left-hand side provided that 
$z\mu\hat\ep^2$
is small enough compared to ${\delta}
\hat\theta$.
Therefore, we eventually obtain \eqref{eq:totenergyestimate0} by substituting
a numerical constant times $\mu\theta\hat\delta\hat\ep^2$ for $\hat{z}$ in the definition of $\zeta$ here above.
 \medskip

The next step consists in expanding  $\int_\Omega F(\fy,\Delta\fy,T)\overline{T}$,  the first term in the right-hand side of \eqref{eq:totenergyestimate0}, and estimating each term 
by a power of $E_1(\overline T,\varphi)$. This is made in the following. 
\begin{lemma}[final a priori estimate]  In addition to the notations in\-tro\-du\-ced in Lemma \ref{lem:E1} (see p.~\pageref{lem:E1}),
let us define the constants $A$, $B$, $C$, $D$, by
\begin{align*}
A:=\,&A^0+\beta  \alpha %
(1+
\mu^{1/2}) W'_\infty+ \hat\delta\alpha^2\mu^{3/2},\\
C:=\,&C^0+\hat\theta \alpha \mu^{1/2} W'_\infty \|\widetilde T\|_{H^1}^2
+\hat\theta \alpha |W'_\infty|^3+ \hat\theta \alpha \ep^2\mu^{1/2} \|\widetilde T\|_{H^1}^4,\\
D:=\,&D^0+\ep^2(\beta\theta)^{1/2}\alpha^2\hat\delta\mu\,(1 
+\alpha^3\hat\delta^3\mu^3 (\beta\theta)^{3/2})\\
&+\ep\beta \mu (\alpha^{3}\hat\delta \theta)^{1/2}(1+\ep^3  \mu^{3} (\alpha^{9}  \theta^{3}\hat\delta^{9})^{1/2})
\end{align*}
Then the first $n$ modes $(T=\overline T+\widetilde{T},\varphi)=(T^n=\overline T^n+\widetilde{T},\varphi^n)$ solutions to~\eqref{galerkin-fy}-\eqref{galerkin-T} satisfy the energy estimate
\begin{multline}
\displaystyle\frac{\dif}{\dif t}E_1^*(\overline T,\varphi)+\tfrac12 {\delta}{\hat\theta} \|\nabla T\|_{L^2}^2+\hat\alpha\|\partial_t\varphi\|_{L^2}^2+\tfrac12
\ep^4 \delta  \hat\mu \hat\theta |\nabla\Delta\varphi\|_{L^2}^2\\
\lesssim AE_1^*(\overline T,\varphi)+D(B+E_1^*(\overline T,\varphi))^3+C,
\label{eq:totenergyestimate}\end{multline}
for all $t\in(0,t_n)$.
\label{lem:FinalEnergyEstimate}\end{lemma}

{\bf Proof:} Recalling the definition of $F$ in \eqref{eq:systsimp1F} on p.~\pageref{eq:systsimp1F}, we have
\begin{equation*}
|F(\fy,\Delta\fy,T)|\leq
 \underbrace{
2
\alpha |W'_\infty|^2+
\alpha\gamma\nu'_\infty W'_\infty |T|}_{f_1}+\underbrace{
2
\alpha\ep^4\,|\Delta\varphi|^2}_{f_2}+\underbrace{
\alpha\ep^2\gamma\nu'_\infty |\Delta\varphi||T|}_{f_3}.
\end{equation*}
\begin{enumerate}
\item 
Since
$$2 |W'_\infty|^2 \overline T \leq |W'_\infty|^3 + W'_\infty \overline T^2\;\mbox{and}\quad
2|T|\overline T \leq 3 \overline T^2+\widetilde{T}^2\,,$$
we readily have
$$
\hat\theta \int_\Omega f_1 \overline T \lesssim 
\hat\theta \alpha 
(1+
\gamma\nu'_\infty) W'_\infty \|\overline T\|_{L^2}^2+
\hat\theta \alpha \gamma\nu'_\infty W'_\infty \|\widetilde T\|_{L^2}^2+
\hat\theta \alpha |W'_\infty|^3\,, 
$$
hence
\begin{equation}
\hat\theta \int_\Omega f_1 \overline T \lesssim 
a_1\, E_1^*(\overline T,\varphi)
+c_1
,\label{estimF0}
\end{equation}
with 
$$a_1:=\beta  \alpha %
(1+
\gamma\nu'_\infty) W'_\infty\,,\quad
c_1:=\hat\theta \alpha \gamma\nu'_\infty W'_\infty \|\widetilde T\|_{H^1}^2
+\hat\theta \alpha |W'_\infty|^3\,.$$

\item By the Cauchy-Schwarz inequality, we have 
$$\int_\Omega |\Delta\varphi|^2 |\overline T| \leq \|\Delta\varphi\|_{L^4}^2\|\overline T\|_{L^2}\,.$$
Using the interpolation inequality  (due to H\"older)
\begin{equation*}
\|v\|_{L^4}\leq \|v\|_{L^2}^{1/4}\|v\|_{L^6}^{3/4},
\end{equation*}
and the Sobolev inequality 
$$\|v\|^2_{L^6}\lesssim \|v\|^2_{L^2}+\|\nabla v\|^2_{L^2}$$
(equivalent to the Sobolev embedding $H^1\hookrightarrow L^6$ already mentioned before),
we thus infer that
\begin{align}
\int_\Omega |\Delta\varphi|^2 |\overline T| &
\lesssim \|\Delta\varphi\|_{L^2}^{2}\|\overline T\|_{L^2}+\|\Delta\varphi\|_{L^2}^{1/2}\|\overline T\|_{L^2}\|\nabla\Delta\varphi\|^{3/2}_{L^2}\label{HIS}.
\end{align}
As a consequence of \eqref{HIS} and Young's inequality (Eq. \eqref{eq:Young}
with a factor $\lambda$ to be determined afterwards, and again $p=4/3$), we get
\begin{align*}
\hat\theta \int_\Omega f_2 \overline T
\lesssim\,&
\hat\theta\alpha \ep^4\,\|\Delta\varphi\|_{L^2}^{2}\|\overline T\|_{L^2}
			+
			\hat\theta\alpha \ep^4\,\left(
			\hat\lambda^3\,\|\Delta\varphi\|_{L^2}^2\|\overline T\|_{L^2}^4+
			\lambda\,\|\nabla\Delta\varphi\|^{2}_{L^2}\right).
\end{align*}
So, using that
\begin{equation}\label{eq:Deltafy}
\|\Delta\fy\|_{L^2}^2\lesssim \hat\ep^2\hat\delta\alpha\mu\theta E_1(\overline{T},\fy)\,,\quad
\|\overline T\|_{L^2}^2\lesssim \beta\theta E_1(\overline{T},\fy)\,,
\end{equation}
we arrive at
\begin{multline}
\hat\theta\int_\Omega f_2 \overline T\lesssim
\alpha^2\ep^2\hat\delta\mu (\beta\theta)^{1/2} E_1(\overline T,\varphi)^{3/2}\\
+
\hat\lambda^3 \alpha^2\ep^2\hat\delta\mu \beta^2\theta^2 E_1(\overline T,\varphi)^3+
\lambda \hat\theta\alpha \ep^4\|\nabla\Delta\varphi\|^{2}_{L^2}
.
\label{estimFf1}\end{multline}
We can now choose $\lambda$ 
for the last term in \eqref{estimFf1} to be absorbed by
$\ep^4\delta\hat\mu\hat\theta\|\nabla\Delta\varphi\|_{L^2}^2$, the last term on the left-hand side of 
\eqref{eq:totenergyestimate0} on p.~\pageref{eq:totenergyestimate0}: it suffices to take 
$\lambda$ small enough compared to $\hat\alpha\delta\hat\mu$.
Hiding the multiplicative constant in the $\lesssim$ sign, we thus get from~\eqref{estimFf1} that
\begin{equation}
\hat\theta \int_\Omega f_2  \overline T- \tfrac12 \ep^4\delta\hat\mu\hat\theta\|\nabla\Delta\varphi\|_{L^2}^2\lesssim 
d_2 E_1^*(\overline T,\varphi)^3
\label{estimF1}\end{equation}
with 
$$d_2:=\ep^2(\beta\theta)^{1/2}\alpha^2\hat\delta\mu\,(1 
+\alpha^3\hat\delta^3\mu^3 (\beta\theta)^{3/2})\,.$$

\item
The way of estimating $
\int_\Omega f_3 \overline T$  is very similar to the  one for $
\int_\Omega f_2 \overline T$. Indeed, we have
$$2 \int_\Omega |\Delta\fy| |T| \overline T\,\leq 
\int_\Omega |\Delta\varphi||\overline T|^2+\int_\Omega |\Delta\varphi||\widetilde T|^2\,,$$
with 
$$\int_\Omega |\Delta\varphi||\widetilde T|^2\,\leq\,\|\Delta\varphi\|_{L^2}\|\widetilde T\|_{H^1}^2$$
by Cauchy-Schwarz and the Sobolev embedding $H^1\hookrightarrow L^4$,
and 
\begin{align*}
\int_\Omega |\Delta\varphi||\overline T|^2\lesssim &
\|\Delta\varphi\|_{L^2}\|\overline T\|_{L^2}^2
+\|\Delta\varphi\|_{L^2}\|\overline T\|_{L^2}^{1/2}\|\nabla \overline T\|_{L^2}^{3/2},
\end{align*}
by exchanging the role of $\Delta\varphi$ and $\overline T$ in~\eqref{HIS}.
Therefore, using again \eqref{eq:Young} we have
\begin{align*} \int_\Omega |\Delta\fy| |T| \overline T\,&\lesssim \|\widetilde T\|_{H^1}^4 + \|\Delta\fy\|_{L^2}^2+\|\Delta\varphi\|_{L^2}\|\overline T\|_{L^2}^2\\
&+\hat\lambda^3 \|\Delta\varphi\|_{L^2}^4\|\overline T\|_{L^2}^{2}+\lambda\,\|\nabla \overline T\|_{L^2}^{2}
\end{align*}
for all positive $\lambda$.
Using again \eqref{eq:Deltafy} we thus infer that
\begin{align}
\hat\theta\int_\Omega f_3 \overline T&\lesssim \hat\theta \alpha \ep^2\mu^{1/2} \|\widetilde T\|_{H^1}^4  + \hat\delta\alpha^2\mu^{3/2} E_1(\overline{T},\fy)  \nonumber\\
& + \ep\beta \mu \alpha^{3/2} (\hat\delta \theta)^{1/2} E_1(\overline{T},\fy)^{3/2} + 
\hat\lambda^3 \hat\ep^2 \beta \alpha^3 \mu^{5/2} \theta^2\hat\delta^2 E_1(\overline{T},\fy)^{3} \nonumber\\
& +\lambda \hat\theta \alpha \ep^2 \mu^{1/2} \|\nabla \overline T\|_{L^2}^{2} \,.
\label{estimFf3}
\end{align} 
Now, if we want to let the last term in \eqref{estimFf3}
 to be absorbed by the left-hand side term $\delta\hat\theta \|\nabla \overline T\|_{L^2}^{2}$  in \eqref{eq:totenergyestimate0}, we choose $\lambda$ small enough compared to
$\delta\hat\alpha\hat\ep^2\hat\mu^{1/2}$.
This yields
\begin{align}
\hat\theta\int_\Omega f_3 \overline T\,-\,\tfrac12 \delta\hat\theta \|\nabla \overline T\|_{L^2}^{2} &\lesssim \hat\theta \alpha \ep^2\mu^{1/2} \|\widetilde T\|_{H^1}^4  + \hat\delta\alpha^2\mu^{3/2} E_1(\overline{T},\fy)  \nonumber\\
& + \ep\beta \mu \alpha^{3/2} (\hat\delta \theta)^{1/2} E_1(\overline{T},\fy)^{3/2}\\ \nonumber
&+ \ep^4 \beta \alpha^6 \mu^{4} \theta^2\hat\delta^5 E_1(\overline{T},\fy)^{3},\nonumber
\end{align}
hence
\begin{equation}
\hat\theta\int_\Omega f_3 \overline T\,-\,\tfrac12 \delta\hat\theta \|\nabla \overline T\|_{L^2}^{2} \lesssim c_3  + a_3 E_1(\overline{T},\fy)  + d_3 E_1^*(\overline{T},\fy)^{3},
\label{estimF3}
\end{equation}
with 
$$a_3:=\hat\delta\alpha^2\mu^{3/2}\,,\; c_3:=\hat\theta \alpha \ep^2\mu^{1/2} \|\widetilde T\|_{H^1}^4\,,$$
and
$$
d_3:=\ep\beta \mu (\alpha^{3}\hat\delta \theta)^{1/2}(1+\ep^3  \mu^{3} (\alpha^{9}  \theta^{3}\hat\delta^{9})^{1/2})\,.$$

\end{enumerate} 
Finally, adding together the estimates in~\eqref{estimF0}-\eqref{estimF1}-\eqref{estimF3} 
and the energy estimate in~\eqref{eq:totenergyestimate0} gives 
\eqref{eq:totenergyestimate} with
$$A:=A^0+a_1+a_3\,,\;C:=C^0+c_1+c_3\,,\;D:=D^0+d_2+d_3\,.$$

\medskip

Recall that $T$ and $\varphi$ in~\eqref{eq:totenergyestimate} are actually $T^n$ and $\varphi^n$, and of course depend on the rank $n$ in the Galerkin approximation. However, their initial data $T_0\in {L^2}$ and $\varphi_0\in {H^2_{\mathrm {n}}}$ are independent of $n$, and so is the initial energy $E_1(\overline T_0,\varphi_0)$. 
By the energy estimate~\eqref{eq:totenergyestimate}, there is a uniform time $t^*>0$ of existence  of $T^n$ and $\varphi^n$ 
-~as solutions of the ODEs  \eqref{eq:ODE}~- such that 
$E(\overline T^n,\varphi^n)$ is bounded in $L^\infty(0,t^*)$, while $\|\nabla\overline T^n\|_{L^2}$, $\|\partial_t\varphi^n\|_{L^2}$ and $\|\nabla\Delta\varphi^n\|_{L^2}$ are bounded in $L^2(0,t^*)$.
Therefore, by the usual compactness theorems,  up to extracting subsequences we have 
\begin{equation*}
\fy_n \rightarrow \fy 
\end{equation*}
in $L^\infty(0,t^*;{H^2_{\mathrm {n}}})$ weak-*, in $L^2(0,t^*;H^3_{\mathrm {n}})$ weak, in $L^2(0,t^*;L^2)$ strong, and also almost everywhere,
\begin{equation*}
\dt \fy_n\rightarrow \dt \fy 
\end{equation*}
in $L^\infty(0,t^*;L^2)$ weak-*,
\begin{equation*}
\overline T_n \rightarrow \overline T
\end{equation*}
in $L^\infty(0,t^*;{L^2})$ weak-*, in $L^2(0,t^*;{H^1_0})$ weak, in $L^2(0,t^*;{L^2})$ strong, and almost everywhere.
Furthermore, since 
$$\dt \fy \in L^2(0,t^*;L^2)\quad \mbox{and}\quad \fy \in L^2(0,t^*;{H^2_{\mathrm {n}}})\,,$$ 
we have by the Aubin-Simon lemma that
\begin{equation*}
\fy\, \in \mathscr C(0,t^*;{H^1}).
\end{equation*}
Moreover,  
the definition of $\overline T^n$ in \eqref{galerkin-T} shows that
 the sequence $(\dt \overline T^n)$ is bounded in the dual of $L^2(0,t^*;H^{1}_0)$, 
which implies that
\begin{equation*}
\overline T\, \in \mathscr C(0,t^*;{L^2}).
\end{equation*}
Passing to the limit in~\eqref{galerkin-fy} and~\eqref{galerkin-T} we conclude that $(T,\varphi)$ is a weak solution to~\eqref{eq:systsimp1} on $[0,t^*)$.
\medskip

\begin{remark}[time of existence] 
By  the energy estimate in~\eqref{eq:totenergyestimate} and an elementary ODE argument we see that the time of existence $t^*$ is bounded from below by $t^*_1:=1/(D E_1(\overline T_0,\varphi_0)^2)$ (up to a multiplicative constant only depending on~$\Omega$ and the dimension~$d$). This yields two comments.
\begin{enumerate}
\item For fixed parameters $\alpha,\beta,\gamma,\delta,\ep,\theta$, by definition of $E_1$, $t^*_1$ tends to infinity when the initial data tend to a constant state $(\overline T_0\equiv 0,\varphi_0\equiv 0)$ or 
$(\overline T_0\equiv 0,\varphi_0\equiv 1)$. In other words, $t^*_1$ is all the more larger as the initial state is close to a single phase.
\item We can also examine how $t^*_1$ varies with respect to the parameters $\alpha$, $\beta$, $\gamma$, $\delta$, $\ep$, $\theta$ for fixed initial data. 
By inspection of all coefficients in the definition of $D$ (p.~\pageref{eq:totenergyestimate}) we see that $D$ is bounded on compact subsets of 
$$\{(\alpha,\beta,\gamma,\delta,\ep,\theta)\,;\;\alpha\geq 0, \beta \geq 0, \gamma\geq 0, \delta>0, \ep>0,\theta\geq 0\}\,,$$
and tends to zero in either one of the limits
$$\beta\to 0 \mbox{ or } \gamma\to 0 \mbox{ or } \theta\to 0 \mbox{ or } \delta\to \infty\,.$$
Note that, by definition of $\gamma$ and $\delta$ (see Eq.~\eqref{eq:newparam} on p.~\pageref{eq:newparam}),
$1/\Stefan=\beta\gamma\theta$, $1/\Peclet=\beta\delta$, so that
$\Stefan$ tends to infinity in either one of the first three limits, whereas
$\Peclet$ can tend to infinity  (if $\beta\to 0$ and $\delta$ fixed), tend to zero
(if $\delta\to \infty$ {and} $\beta$ fixed), or can be kept bounded and bounded away from zero (if 
$\beta\to 0$ \emph{and} $\delta \to \infty$ with $\beta$ and $1/\delta$ of the same order).
Unfortunately, it turns out that neither one of these limits is compatible with keeping $E_1(\overline T_0,\varphi_0)$ bounded (which would imply $t^*\to \infty$), or at least
$D E_1(\overline T_0,\varphi_0)^2$ bounded  (which would give a uniform lower bound for $t^*$), independently of $\|\overline T_0\|_{L^2}$ and $ \|\Delta \fy_0\|_{L^2}$.
Indeed, we have
$$2\beta \theta E_1(\overline T_0,\varphi_0)\geq \|\overline T_0\|_{L^2}^2\;,\quad 
2 \alpha\,\gamma^2 |\nu'_\infty|^2 \theta E_1(\overline T_0,\varphi_0)\geq  \ep^{2} \delta \|\Delta \fy_0\|_{L^2}^2,$$
and looking closer at $D$ we have
$$D\stackrel{\beta\to 0}{=}O(\beta^{1/2})\,,\;
D\stackrel{\gamma\to 0}{=}O(\gamma^{2})\,,\;
D\stackrel{\hat\delta\to 0}{=}O(\hat\delta^{1/2})\,,\;
D\stackrel{\theta\to 0}{=}O(\theta^{1/2})\,.
$$
So, the limit $\beta\to 0$ penalizes $\|\overline T_0\|_{L^2}$, 
while $\theta \to 0$ penalizes both $\|\overline T_0\|_{L^2}$ and $\|\Delta \fy_0\|_{L^2}$,
and both $\gamma\to 0$ and $\delta\to \infty$ penalize $\|\Delta \fy_0\|_{L^2}$.
Of course, if we restrict our initial data to $\overline T_0\equiv 0$, the limit $\beta\to 0$  is allowed and thus gives 
$t^* \to \infty$.
This is no surprise because in this case we are basically left with the  Allen-Cahn equation
 (the first equation in \eqref{eq:systsimpbis} with $T= \widetilde{T}$), 
for which global existence is well-known. 
\end{enumerate}

\end{remark}

\subsection{Continuous dependence on the data}

This final section is devoted to the proof of Theorem~\ref{theorem-unicitesolution}, which gives both uniqueness and continuous dependence of the initial data.  The tools are basically
the same as in the  existence proof (Theorem~\ref{theorem-localweaksolution}), namely,
energy estimates and various inequalities of  Sobolev and/or interpolation type. 
However, the details are longer and somehow more technical. To gain some simplicity in the exposure we set all six parameters $\alpha$, $\beta$, $\gamma$, $\delta$, $\ep$, and $\theta$ equal to $1$ in \eqref{eq:systsimp1}, and 
the $L^\infty$ norms of derivatives of $W$ and $\nu$ will be systematically hidden in the $\lesssim$ sign.
\medskip

{\bf Step 1.}
We assume that $(T_i,\varphi_i)$, $i\in\{1,2\}$ are two weak solutions to~\eqref{eq:systsimp1} on $[0,t^*)$ for some positive $t^*$. We are going to use repeatedly the notation
\begin{equation*}
[G]=G_2-G_1
\end{equation*}
for the difference of any two quantities $G_2$ and $G_1$. The differences $[T]$ and $[\varphi]$ then satisfy the system
\begin{equation}\label{eq:systsimp1diff}
\left\{
\begin{aligned}
& \partial_t[\varphi] -\Delta [\varphi]+[W'(\varphi)]=[\nu'(\varphi)\,T], \\
& \partial_t [T] +[\nu'(\varphi)\partial_t\varphi]-\Delta [T]=[F],\\
&F_i:=(\Delta\varphi_i-W'(\varphi_i))^2+(\Delta\varphi_i-W'(\varphi_i))\nu'(\varphi_i)T_i.
\end{aligned}
\right.
\end{equation}
By analogy with $E_0$ defined in \eqref{eq:E0} we set 
\begin{equation*}
\displaystyle e_0([T],[\varphi])=\int_\Omega \tfrac12 
[T]^2+\tfrac12 
|\nabla[\varphi]|^2+ [\varphi][W'(\varphi)].
\end{equation*}
Similarly to the energy estimate~\eqref{eq:energyCaginalpcoeff} for $E_0$ (except that the potential energy has been discarded here), we derive the identity
\begin{multline}
\displaystyle\frac{\dif}{\dif t}e_0([T],[\varphi]) + \int_\Omega |\nabla[T]|^2+\int_\Omega|\partial_t[\varphi]|^2\\
=\int_\Omega [\nu'(\varphi)T][\partial_t\varphi]-[\nu'(\varphi)\partial_t\varphi] [T]
+\int_\Omega[\partial_t\varphi][W'(\varphi)]+\int_\Omega [F] [T].
\label{eq:e01}\end{multline}
Using the formula of differentiation
\begin{equation*}
[GH]=\<G\> [H]+[G]\<H\>,\quad \<G\>:=\frac{G_2+G_1}{2},
\end{equation*}
and the mean value theorem (for $\nu'$)
we obtain
\begin{align}
[\nu'(\varphi)T][\partial_t\varphi]-[\nu'(\varphi)\partial_t\varphi] [T]
=&[\nu'(\varphi)](\<\partial_t\varphi\>[T]-\<T\>[\partial_t\varphi])\nonumber\\
\lesssim &\,
|\<\partial_t\varphi\>||[\varphi]||[T]|+
|\<T\>||[\varphi]||[\partial_t\varphi]|.
\label{eq:e02}\end{align}
In order to estimate the first term of the right-hand side in~\eqref{eq:e02}, we can use the following inequality,
whose proof is postponed to the appendix,
\begin{equation}
\int_\Omega a b c
\leq \ep \|\nabla b\|_{L^2}^2+\ep \|\nabla c\|_{L^2}^2+C_\ep(1+ \|a\|_{L^2}^4)(\|b\|_{L^2}^2+\|c\|_{L^2}^2),
\label{ineq:abc}\end{equation}
satisfied for all $a\in L^2(\Omega)$, $b,c\in H^1(\Omega)$, $\ep>0$. Here above and in what follows, $C_\ep$ denotes a constant depending only on $\ep$, $\Omega$ %
and $d$. The inequality~\eqref{ineq:abc} gives
\begin{align}
\int_\Omega 
|\<\partial_t\varphi\>||[\varphi]||[T]|&
\leq 
\ep \|\nabla[T]\|_{L^2}^2+\ep\|\nabla[\varphi]\|_{L^2}^2
\nonumber\\ 
&
+C_\ep(\|\<\partial_t\varphi\>\|_{L^2}^4+1)(\|[T]\|_{L^2}^2+\|\nabla[\varphi]\|_{L^2}^2).
\label{eq:e03}\end{align}
To obtain a bound on the second term in the right-hand side of~\eqref{eq:e02}, we are going to use Agmon's inequality \eqref{eq:Agmon}
under the form
\begin{equation}
\|b\|_{L^\infty}^2\lesssim \|b\|_{L^2}^2+\|\nabla b\|_{L^2}\|\Delta b\|_{L^2},
\label{eq:AgmonBis}\end{equation}
which implies (by Cauchy-Schwarz) that
\begin{equation}
\int_\Omega |abc|\leq \ep\|c\|_{L^2}^2+C_\ep \|a\|_{L^2}^2(\|b\|_{L^2}^2+\|\nabla b\|_{L^2}^2+\|\Delta b\|_{L^2}^2).
\label{ineq:a0b0c0}\end{equation}
We thus find that
\begin{align}
\int_\Omega |\<T\>||[\varphi]||[\partial_t\varphi]|
&\leq \ep \|[\partial_t\varphi]\|_{L^2}^2\nonumber \\
&+C_\ep  \|\<T\>|^2_{L^2}(\|[\varphi]\|_{L^2}^2+\|\nabla [\varphi]\|_{L^2}^2+\|\Delta [\varphi]\|_{L^2}^2)\label{ineq:a0b0c00}.
\end{align}
Together with~\eqref{eq:e03},~\eqref{ineq:a0b0c00} gives
\begin{multline}
\int_\Omega [\nu'(\varphi)T][\partial_t\varphi]-[\nu'(\varphi)\partial_t\varphi] [T] \leq  \ep \|\nabla[T]\|_{L^2}^2+\ep \|[\partial_t\varphi]\|_{L^2}^2\\
+C_\ep(1+\|\<\partial_t\varphi\>\|_{L^2}^4+\|\<T\>|^2_{L^2})(\|[T]\|_{L^2}^2+\|[\varphi]\|_{L^2}^2+\|\nabla [\varphi]\|_{L^2}^2+\|\Delta [\varphi]\|_{L^2}^2).
\label{eq:e04}\end{multline}
Let us introduce the shorthand $\overline{z}_0:=1+\|\<\partial_t\varphi\>\|_{L^2}^4+ \|\<T\>\|^2_{L^2}$, and 
\begin{equation*}
\displaystyle e_1([T],[\varphi])=e_0([T],[\varphi])+\mu \|\Delta[\varphi]\|_{L^2}^2\,,\quad
e_2([T],[\varphi])=e_1([T],[\varphi])+\|[\varphi]\|_{L^2}^2,
\end{equation*}
where $\mu$ is a parameter that will be chosen later. 
The final estimate we are aiming at reads
\begin{equation}
\displaystyle\frac{\dif}{\dif t}e_2\leq \overline{z} e_2,
\quad e_2=e_2([T],[\varphi]),
\label{eq:e12}\end{equation}
for some $\overline{z}\in L^1(0,t^*)$ depending
continuously
on the norms of $T_i$ in $L^\infty(0,t^*;H^1
)$, $L^2(0,t^*;L^2
)$ and the norms of $\varphi_i$ in $L^\infty(0,t^*;H^2
)$, $L^2(0,t^*;H^3
)$, $i\in\{1,2\}$\footnote{In particular, using the equation for $\varphi$ in~\eqref{eq:systsimp1}, $\overline{z}$ may depend on the norm of $\partial_t\varphi$ in $L^\infty(0,t^*;L^2
)$}. 
Once we have~\eqref{eq:e12}, the continuous dependence on the data as expressed in~\eqref{eq:continuousdependencedata} will be clear because, by straightforward integration
\begin{equation*}
e_2(t)\leq e^{\int_0^t \overline{z}(s)ds} e_2(0),
\end{equation*}
and by definition 
$e_2$ controls $\|[\fy|\|_{H^2}$ and $\|[T]\|_{L^2}$.
\medskip

The first step in the derivation of \eqref{eq:e12} follows from
the identity in~\eqref{eq:e01} and the inequality in~\eqref{eq:e04}, which together give 
\begin{multline}
\displaystyle\frac{\dif}{\dif t}e_0([T],[\varphi]) +
\|\nabla[T]\|_{L^2}^2+
\|\partial_t[\varphi]\|_{L^2}^2\leq \ep \|\nabla[T]\|_{L^2}^2+\ep \|[\partial_t\varphi]\|_{L^2}^2\\
+C_{\ep,\mu}\overline{z}_0 e_1([T],[\varphi])
+\int_\Omega[\partial_t \varphi][W'(\varphi)]+\int_\Omega [F] [T],
\label{eq:e11}\end{multline}
where $\overline{z}_0\in L^1(0,t^*)$, 
and $C_{\ep,\mu}:=\max(1,1/\mu)\,C_\ep$.
We will prove~\eqref{eq:e12} by using an additional energy estimate for $\|\Delta[\varphi]\|_{L^2}^2$ and some estimates for the remaining terms
\begin{equation*}
\int_\Omega[\partial_t \varphi][W'(\varphi)]
\quad \text{and} \quad
\int_\Omega [F] [T]
\end{equation*}
in the right hand-side of~\eqref{eq:e11}. 
\medskip

{\bf Step 2.}
By \eqref{eq:elem}, 
\begin{equation}
\int_\Omega [\partial_t \varphi][W'(\varphi)]\leq  \ep \|[\partial_t\varphi]\|_{L^2}^2+\overline z_1  \|[\varphi]\|_{L^2}^2.
\label{eq:e14}\end{equation}

{\bf Step 3.}
There are several factors in the expansion of the term $\int_\Omega [F][T]$ that we bound successively and, to some extent, by decreasing order of difficulty.

By using the expansion  formula
\begin{equation*}
[GHI]=\<G\>\,\<H\>[I]+\<G\>\,\<I\>[H]+\<HI\>[G],
\end{equation*}
we obtain
\begin{equation}
[\Delta\varphi \nu'(\varphi) T][T]=\<\Delta\varphi\>\,\<\nu'(\varphi)\>[T]^2+\<\Delta\varphi\>\,\<T\>[\nu'(\varphi)][T]+\<\nu'(\varphi)T\>[\Delta\varphi][T].
\label{eq:e14bis}\end{equation}
By the inequality in~\eqref{ineq:abc} 
we have
\begin{equation}
\int_\Omega \<\Delta\varphi\>\,\<\nu'(\varphi)\>[T]^2\leq\ep\|\nabla[T]\|_{L^2}^2+C_\ep(1+\|\<\Delta\varphi\>\|_{L^2}^4)\|[T]\|_{L^2}^2.
\label{eq:e15}\end{equation}
and similarly
\begin{align}
\int_\Omega \<\nu'(\varphi)T\>[\Delta\varphi][T]&\leq\ep\|\nabla[T]\|_{L^2}^2+\ep\|\Delta[\varphi]\|_{L^2}^2\nonumber \\
&+C_\ep(1+\|\<T\>\|_{L^2}^4)(|[T]\|_{L^2}^2+\|[\varphi]\|_{L^2}^2).
\label{eq:e160}\end{align}
To estimate the remaining term in~\eqref{eq:e14bis}, we apply the inequality (proved in the appendix)
\begin{align}
\int_\Omega |abcd|\leq & \ep\|\nabla d\|_{L^2}^2+\ep\|\Delta c\|_{L^2}^{2} \nonumber \\
& + C_\ep(1+\|a\|_{L^2}^4+\|b\|_{L^2}^{4}+\|b\|_{L^2}^{2}\|\nabla b\|_{L^2}^{2})(\|\nabla c\|_{L^2}^{2}+\|d\|_{L^2}^2)
\label{ineq:abcd}\end{align}
to $a=\<\Delta\varphi\>$, $b=\<T\>$, $c=[\varphi]$, $d=[T]$. This gives
\begin{multline*}
\int_\Omega |\<\Delta\varphi\>\,\<T\>[\nu'(\varphi)][T]|
\leq \ep\|\nabla [T]\|_{L^2}^2+\ep\|\Delta [\varphi]\|_{L^2}^{2}\\
+C_\ep(1+\|\<\Delta\varphi\>\|_{L^2}^4+\|\<T\>\|_{L^2}^{4}+\|\<T\>\|_{L^2}^{2}\|\nabla \<T\>\|_{L^2}^{2})(\|\nabla [\varphi]\|_{L^2}^{2}+\|[T]\|_{L^2}^2).
\end{multline*}
With~\eqref{eq:e15} and~\eqref{eq:e160}, 
we arrive at
\begin{equation}
\int_\Omega |[\Delta\varphi \nu'(\varphi) T][T]|\leq \ep\|\nabla [T]\|_{L^2}^2+\overline{z}_{2}
(\|[\varphi]\|_{L^2}^2+e_1).
\label{eq:e16}\end{equation}
Replacing $\Delta\varphi$ by $W'(\varphi)$ in~\eqref{eq:e16} and using the estimate $|[W'(\varphi)]|\leq W''_\infty|[\varphi]|$, we will obtain a similar estimate for the quadrilinear term $|[W'(\varphi) \nu'(\varphi) T][T]$. For the trilinear term 
\begin{equation*}
[(\Delta\varphi-W'(\varphi))^2][T]=[|\Delta\varphi|^2][T]+[W'(\varphi)|^2][T]-2[\Delta\varphi W'(\varphi)][T],
\end{equation*}
by proceeding as in Step 2, we also have an estimate as in~\eqref{eq:e16}. We thus conclude with Step 2 that
\begin{equation}
\int_\Omega |[\varphi][\partial_t W'(\varphi)]|+\int_\Omega |[F] [T]|\leq \ep\|\nabla [T]\|_{L^2}^2+\overline{z}_3
(\|[\varphi]\|_{L^2}^2+e_1).
\label{eq:e17}\end{equation}
Inserting~\eqref{eq:e17} in~\eqref{eq:e11} gives
\begin{equation}
\displaystyle\frac{\dif}{\dif t}e_0+ \|\nabla[T]\|_{L^2}^2+\|\partial_t[\varphi]\|_{L^2}^2\\
\leq \ep \|\nabla[T]\|_{L^2}^2+\ep \|[\partial_t\varphi]\|_{L^2}^2+\overline{z}_4
(\|[\varphi]\|_{L^2}^2+e_1).
\label{eq:e18}\end{equation}
To eventually obtain~\eqref{eq:e12}, we still have to estimate $\displaystyle\frac{\dif}{\dif t}\|\Delta[\varphi]\|_{L^2}^2$.
\medskip

{\bf Step 4.}
Taking $-\Delta^2[\varphi]$ as a test function\footnote{To be exact, we take $-\Delta^2[\varphi]^n$ as a test function, where $[\varphi]^n:=\sum_{i=1}^n ([\varphi],\varphi_i)\varphi_i$ is the Galerkin approximation of rank $n$ to $[\varphi]$ and, in a final step, pass to the limit $n\to+\infty$ in the estimates obtained for $[\varphi]^n$.} in the equation satisfied by $[\varphi]$ in~\eqref{eq:systsimp1diff}, we obtain
\begin{equation}
\frac{\dif\;}{\dif t}\|\Delta[\varphi]\|_{L^2}^2+\|\nabla\Delta[\varphi]\|_{L^2}^2\leq \int_\Omega |[\nabla(\nu'(\varphi)T)][\nabla\Delta\varphi]|
+\int_\Omega |[\nabla W'(\varphi)][\nabla\Delta\varphi]|.
\label{eq:ephi1}\end{equation}
We have
\begin{equation}
\int_\Omega |[\nabla W'(\varphi)][\nabla\Delta\varphi]|\leq \ep \|[\nabla\Delta\varphi]\|_{L^2}^2+C_\ep \|[\varphi]\|_{L^2}^2.
\label{eq:ephi2}\end{equation}
Besides, we have the expansion
\begin{equation}
|[\nabla(\nu'(\varphi)T)]=\<\nu'(\varphi)\>[\nabla T]+\<\nabla T\>[\nu'(\varphi)]+\<\nabla\nu'(\varphi)\>[T]+\<T\>[\nabla\nu'(\varphi)]
\label{eq:ephi3}\end{equation}
which has to be multiplied by $[\nabla\Delta\varphi]$.
We have
\begin{equation}
\int_\Omega |\<\nu'(\varphi)\>[\nabla T][\nabla\Delta\varphi]|\leq\ep \|[\nabla\Delta\varphi]\|_{L^2}^2+C_\ep \|[T]\|_{L^2}^2
\label{eq:ephi4}\end{equation}
for the first term, and by~\eqref{ineq:a0b0c0}
\begin{align}
\int_\Omega |\<\nabla T\>[\nu'(\varphi)][\nabla\Delta\varphi]|&\leq\ep \|[\nabla\Delta\varphi]\|_{L^2}^2\nonumber\\
&+C_\ep \|\<\nabla T\>\|_{L^2}^2(\|[\varphi]\|_{L^2}^2+\|\nabla[\varphi]\|_{L^2}^2+\|\Delta[\varphi]\|_{L^2}^2)
\label{eq:ephi5}\end{align}
for the second term. To estimate the third term, we apply the inequality (proved in the appendix)
\begin{equation}
\int_\Omega |abc|\leq  \ep\|c\|_{L^2}^2+C_\ep\|\nabla b\|_{L^2}^2+C_\ep (1+\|\nabla a\|_{L^2}^4+\|a\|_{L^2}^4)\|b\|_{L^2}^2
\label{ineq:abc1}\end{equation}
to $a=\<\nabla\varphi\>$, $b=[T]$, $c=[\nabla\Delta\varphi]$. We  obtain
\begin{align}
\int_\Omega|\<\nabla\nu'(\varphi)\>[T][\nabla\Delta\varphi]|&\leq\ep\|[\nabla\Delta\varphi]\|_{L^2}^2+C_\ep \|\nabla [T]\|_{L^2}^2\nonumber\\
+C_\ep (1+\<\|\nabla\varphi\|_{L^2}^4\>+\<\|\Delta\varphi\|_{L^2}^4\>)|[T]\|_{L^2}^2.
\label{eq:ephi6}\end{align}
Since $[\nabla\nu'(\varphi)]=\<\nu''(\varphi)\>[\nabla\varphi]+\<\nabla\varphi\>[\nu'(\varphi)]$, the fourth term 
can be split again. Similarly as in~\eqref{eq:ephi5}, we have
\begin{align}
\int_\Omega |\<T\>\<\nu''(\varphi)\>[\nabla\varphi][\nabla\Delta\varphi]|&\leq\ep \|[\nabla\Delta\varphi]\|_{L^2}^2\nonumber\\
&+C_\ep \|\<T\>\|_{L^2}^2(\|[\varphi]\|_{L^2}^2+\|\nabla[\varphi]\|_{L^2}^2+\|\Delta[\varphi]\|_{L^2}^2).
\label{eq:ephi7}\end{align}
There remains to bound
\begin{equation*}
\int_\Omega |\<T\>\,\<\nabla\varphi\>[\varphi][\nabla\Delta\varphi]|.
\end{equation*}
To this purpose, we apply the inequality (proved in the appendix)
\begin{align}
\int_\Omega |abcd|&\leq  \ep\|d\|_{L^2}^2\nonumber\\
&+C_\ep(\|a\|_{L^2}+\|\nabla a\|_{L^2})(\|b\|_{L^2}^2+\|\nabla b\|_{L^2}^2) (\|c\|_{L^2}^2+\|\nabla c\|_{L^2}^2+\|\Delta c\|_{L^2}^2)
\label{ineq:abcd1}\end{align}
to $a=\<T\>$, $b=\<\nabla\varphi\>$, $c=[\varphi]$, $d=[\nabla\Delta\varphi]$. We then deduce from~\eqref{eq:ephi1}, \eqref{eq:ephi2}, \eqref{eq:ephi4}, \eqref{eq:ephi5}, \eqref{eq:ephi6}, \eqref{eq:ephi7}, taking~$\ep$ small enough, the estimate
\begin{equation}
\displaystyle\frac{\dif }{\dif t}\|\Delta[\varphi]\|_{L^2}^2+\tfrac{1}{2}\|\nabla\Delta[\varphi]\|_{L^2}^2\leq C_1\|\nabla T\|_{L^2}^2+C\overline{z_5}
(\|[\varphi]\|_{L^2}^2+e_1).
\label{eq:ephi8}\end{equation}
Choosing $\mu$ so that $\mu C_1$ is small enough, it then follows from~\eqref{eq:e18} and~\eqref{eq:ephi8} that
\begin{equation}
\displaystyle\frac{\dif}{\dif t}e_1+ \|\nabla[T]\|_{L^2}^2+\|\partial_t[\varphi]\|_{L^2}^2\\
\leq \ep \|\nabla[T]\|_{L^2}^2+\ep \|[\partial_t\varphi]\|_{L^2}^2+\overline{z}_6
(\|[\varphi]\|_{L^2}^2+e_1).
\label{eq:ephi9}\end{equation}
To conclude, we perform an energy estimate on $[\varphi]$. Multiplying the first equation in~\eqref{eq:systsimp1diff} 
by $[\varphi]$  we obtain
\begin{align*}
\frac{\dif\;}{\dif t}\|[\varphi]\|_{L^2}^2+\|\nabla[\varphi]\|_{L^2}^2&\leq C\left(\|[\varphi]\|_{L^2}^2+\|[T]\|_{L^2}^2)+\int_\Omega |\<T\>[\varphi]_{L^2}^2|\right)\\
&\leq C\left(|[\varphi]\|_{L^2}^2+\|[T]\|_{L^2}^2)+\|\<T\>\|_{L^2}\|[\varphi]\|_4^2\right)\\
&\leq C(1+\|\<T\>\|_{L^2}^2)(\|[\varphi]\|_{L^2}^2+\|[T]\|_{L^2}^2)+\|\nabla[\varphi]\|_{L^2}^2).
\end{align*}
Adding with~\eqref{eq:ephi9} we get
\begin{equation*}
\displaystyle\frac{\dif}{\dif t}\left(e_1+\|[\varphi]\|_{L^2}^2\right)\leq \overline{z}_7
(e_1+\|[\varphi]\|_{L^2}^2),
\end{equation*}
which is the claimed estimate \eqref{eq:e12}.
This completes the proof of Theorem~\ref{theorem-unicitesolution}.

\section*{Appendix}
We give here the proof of various functional inequalities used previously.
\medskip

{\bf Proof:}[Proof of~\eqref{ineq:abc}]
$$\int_\Omega a b c
\leq \ep \|\nabla b\|_{L^2}^2+\ep \|\nabla c\|_{L^2}^2+C_\ep(1+ \|a\|_{L^2}^4)(\|b\|_{L^2}^2+\|c\|_{L^2}^2)
$$
for $\Omega\subset \R^d$, $d\leq 3$.
By H\"older's inequality, 
\begin{equation*}
\int_\Omega a b c\leq \|a\|_{L^2}\|b\|_{L^6}\|c\|_{L^3}.
\end{equation*}
In addition, we have the interpolation inequality $\|c\|_{L^3}\leq\|c\|_{L^2}^{1/2}\|c\|_{L^6}^{1/2}$, and applying twice  Young's inequality \eqref{eq:Young} (once with $p=4/3$, once with $p=3/2$), we get
\begin{equation*}
\|b\|_{L^6}\cdot \|c\|_{L^6}^{1/2}\cdot \|a\|_{L^2}\|c\|_{L^2}^{1/2}\leq \ep\|b\|_{L^6}^2+\ep\|c\|_{L^6}^{2}+C_\ep\|a\|_{L^2}^4\|c\|_{L^2}^{2}.
\end{equation*}
We can conclude thanks to the Sobolev inequality $\|b\|_{L^6}^2\lesssim \|b\|_{L^2}^2+\|\nabla b\|_{L^2}^2$.

{\bf Proof:}[Proof of~\eqref{ineq:abcd}]
\begin{align}
\int_\Omega |abcd|\leq & \ep\|\nabla d\|_{L^2}^2+\ep\|\Delta c\|_{L^2}^{2} \nonumber \\
& + C_\ep(1+\|a\|_{L^2}^4+\|b\|_{L^2}^{4}+\|b\|_{L^2}^{2}\|\nabla b\|_{L^2}^{2})(\|\nabla c\|_{L^2}^{2}+\|d\|_{L^2}^2)\nonumber
\end{align}
for $\Omega\subset \R^d$, $d\leq 3$.
By H\"older's inequality, we have
\begin{equation*}
\int_\Omega |a b c d|\leq \|a\|_{L^2}\|b\|_{L^3}\|c\|_{L^\infty}\|d\|_{L^6},
\end{equation*}
where we can apply Agmon's inequality~\eqref{eq:AgmonBis} to $c$:
\begin{equation*}
\|c\|_{L^\infty}\leq C(\|c\|_{L^2}+\|\nabla c\|_{L^2}^{1/2}\|\Delta c\|_{L^2}^{1/2}), 
\end{equation*}
and the interpolation inequality $\|b\|_{L^3}\leq \|b\|_{L^2}^{1/2}\|b\|_{L^6}^{1/2}$. This gives
\begin{equation*}
\int_\Omega |abcd|\leq C(\|c\|_{L^2}+\|\nabla c\|_{L^2}^{1/2}\|\Delta c\|_{L^2}^{1/2})\|a\|_{L^2}\|b\|_{L^2}^{1/2}\|b\|_{L^6}^{1/2}\|d\|_{L^6}.
\end{equation*}
By the Sobolev inequality $\|b\|_{L^6}\leq C(\|b\|_{L^2}+\|\nabla b\|_{L^2})$, it follows that
\begin{align}
\int_\Omega |abcd|\leq &\; C(\|c\|_{L^2}+\|\nabla c\|_{L^2}^{1/2}\|\Delta c\|_{L^2}^{1/2}) \nonumber\\
& \times \|a\|_{L^2}\|b\|_{L^2}^{1/2}(\|b\|_{L^2}^{1/2}+\|\nabla b\|_{L^2}^{1/2})(\|d\|_{L^2}+\|\nabla d\|_{L^2}).
\label{ineq:abcd0}\end{align}
The term with the highest order 
derivatives in~\eqref{ineq:abcd0} is
\begin{align*}
& \|\nabla c\|_{L^2}^{1/2}\|\Delta c\|_{L^2}^{1/2}\|a\|_{L^2}\|b\|_{L^2}^{1/2}\|\nabla b\|_{L^2}^{1/2}\|\nabla d\|_{L^2}
=\\
& \|\nabla d\|_{L^2}\cdot \|\Delta c\|_{L^2}^{1/2}\cdot \|\nabla c\|_{L^2}^{1/2}\|a\|_{L^2}\|b\|_{L^2}^{1/2}\|\nabla b\|_{L^2}^{1/2}\\
&\leq\ep\|\nabla d\|_{L^2}^2+\ep\|\Delta c\|_{L^2}^{2}+C_\ep  \|\nabla c\|_{L^2}^{2}\|a\|_{L^2}^4\|b\|_{L^2}^{2}\|\nabla b\|_{L^2}^{2}
\end{align*}
by Young's inequality in the form already used above, namely $$\alpha_1\alpha_2\alpha_3\leq\ep\alpha_1^2+\ep\alpha_2^4+C_\ep\alpha_3^4\,.$$ We finally obtain~\eqref{ineq:abcd} by using similar estimates for the lower order terms in~\eqref{ineq:abcd0}.

{\bf Proof:}[Proof of~\eqref{ineq:abc1}] 
$$\int_\Omega |abc|\leq  \ep\|c\|_{L^2}^2+C_\ep\|\nabla b\|_{L^2}^2+C_\ep (1+\|\nabla a\|_{L^2}^4+\|a\|_{L^2}^4)\|b\|_{L^2}^2
$$
for $\Omega\subset \R^d$, $d\leq 3$.
As above, applying successively the H\"older, Young, interpolation, and Sobolev inequalities give
\begin{align*}
\int_\Omega |abc| &\leq \|c\|_{L^2}\|b\|_{L^3}\|a\|_{L^6}\\
&\leq \ep \|c\|_{L^2}^2+C_\ep \|b\|_{L^3}^2\|a\|_{L^6}^2\\
&\leq \ep \|c\|_{L^2}^2+C_\ep\|b\|_{L^2}\|b\|_{L^6}\|a\|_{L^6}^2\\
&\leq \ep\|c\|_{L^2}^2+C_\ep\|b\|_{L^6}^2+C_\ep \|a\|_{L^6}^4\|b\|_{L^2}^2\\
&\leq \ep\|c\|_{L^2}^2+ C_\ep\|\nabla b\|_{L^2}^2+C_\ep (1+\|\nabla a\|_{L^2}^4+\|a\|_{L^2}^4)\|b\|_{L^2}^2.
\end{align*}

{\bf Proof:}[Proof of~\eqref{ineq:abcd1}] 
\begin{align*}
\int_\Omega |abcd|&\leq  \ep\|d\|_{L^2}^2\\
&+C_\ep(\|a\|_{L^2}+\|\nabla a\|_{L^2})(\|b\|_{L^2}^2+\|\nabla b\|_{L^2}^2) (\|c\|_{L^2}^2+\|\nabla c\|_{L^2}^2+\|\Delta c\|_{L^2}^2)
\end{align*}
for $\Omega\subset \R^d$, $d\leq 3$.
Again, by the  H\"older, Young, interpolation, Sobolev, and Agmon inequalities, 
\begin{align*}
\int_\Omega |abcd| &\leq \|d\|_{L^2}\|a\|_{L^3}\|b\|_{L^6}\|c\|_{L^\infty}\\
&\leq \ep\|d\|_{L^2}^2+C_\ep\|a\|_{L^3}^2\|b\|_{L^6}^2 \|c\|_{L^\infty}^2\\
\leq  \ep\|d\|_{L^2}^2+& C_\ep(\|a\|_{L^2}+\|\nabla a\|_{L^2})(\|b\|_{L^2}^2+\|\nabla b\|_{L^2}^2)\\
 & \times  (\|c\|_{L^2}^2+\|\nabla c\|_{L^2}^2+\|\Delta c\|_{L^2}^2).
\end{align*}

\section*{Note added in proof}
After completing the manuscript, we became aware of a recent work by Feireisl, Petzeltov\'a and Rocca \cite{FeireislPetzeltovaRocca}, in which they analyze a closely related phase transition model with microscopic movements. 
In particular, they take into account a quadratic term $|\partial_t\varphi|^2$ ($\varphi$ being the order parameter) in the equation for the temperature, like the second equation in our system \eqref{eq:systsimp}, which quite equivalently involves the quadratic term $\Delta\varphi\partial_t\varphi$. They do not address the sharp interface limit, but they study local existence, and even \emph{global existence} of a  certain type of weak solutions, obtained by passing to the limit in a
regularized system, in which the diffusion term $-\Delta T$ is replaced by $-\div((1+\varepsilon T^m)\nabla T)$, $m$ large enough. Global existence of classical solutions, and even of weak solutions according to the notion we have used in the present paper (see Definition~\ref{definition-weaksolution}) remains an open question up to our knowledge.

\section*{Acknowledgments}
The authors warmly thank Giulio Schimperna for drawing their attention to the related work \cite{FeireislPetzeltovaRocca}.
The first and second authors have been partly supported by ANR project 08-SYSC-010
MANIPHYC.


\begin{thebibliography}{99}

\bibitem{AndersonMcFaddenWheeler} (MR1609626)
D.~M. Anderson, G.~B. McFadden, and A.~A. Wheeler,
\newblock \emph{Diffuse-interface methods in fluid mechanics},
\newblock Annu. Rev. Fluid Mech., \textbf{30} (1998), 139--165.

\bibitem{BonettiColliFabrizioGilardi09} (MR2507957)
E.~Bonetti, P.~Colli, M.~Fabrizio, and G.~Gilardi,
\newblock \emph{Existence and boundedness of solutions for a singular phase field
  system},
\newblock {J. Differential Equations}, textbf{246} (2009), 3260--3295.

\bibitem{Caginalp} (MR816623)
G.~Caginalp,
\newblock \emph{An analysis of a phase field model of a free boundary},
\newblock {Arch. Rational Mech. Anal.}, \textbf{92}, (1986), 205--245.

\bibitem{CaginalpChen} (MR1643668)
G.~Caginalp and X.~Chen,
\newblock \emph{Convergence of the phase field model to its sharp interface limits},
\newblock {European J. Appl. Math.}, \textbf{9}, (1998), 417--445.

\bibitem{CavaterraGalGrasselliMiranville10} (MR2577805)
C.~Cavaterra, C.~G. Gal, M.~Grasselli, and A.~Miranville,
\newblock \emph{Phase-field systems with nonlinear coupling and dynamic boundary
  conditions},
\newblock {Nonlinear Anal.}, \textbf{72}, (2010), 2375--2399.

\bibitem{ColliHilhorstIssardRochSchimperna09} (MR2525167)
P.~Colli, D.~Hilhorst, F.~Issard-Roch, and G.~Schimperna,
\newblock \emph{Long time convergence for a class of variational phase-field models},
\newblock {Discrete Contin. Dyn. Syst.}, \textbf{25}, (2009), 63--81.

\bibitem{FeireislPetzeltovaRocca} (MR2535695)
E.~Feireisl,  H.~Petzeltov‡, E.~Rocca,
\newblock \emph{Existence of solutions to a phase transition model with microscopic
 movements},
\newblock{Math. Methods Appl. Sci.},  \textbf{32},  (2009),  no. 11, 1345--1369.

\bibitem{GilbargTrudinger} (MR1814364)
D.~Gilbarg and N.~S. Trudinger,
\newblock ``Elliptic partial differential equations of second order",
\newblock Classics in Mathematics. Springer-Verlag, Berlin, 2001.
\newblock Reprint of the 1998 edition.

\bibitem{GrasselliMiranvilleSchimperna10} (MR2629473)
M.~Grasselli, A.~Miranville, and G.~Schimperna,
\newblock \emph{The {C}aginalp phase-field system with coupled dynamic boundary
  conditions and singular potentials},
\newblock {Discrete Contin. Dyn. Syst.}, \textbf{28}, (2010), 67--98.

\bibitem{GrasselliPetzeltovaSchimperna} (MR2216881)
M.~Grasselli, H.~Petzeltov{\'a}, and G.~Schimperna,
\newblock \emph{Long time behavior of solutions to the {C}aginalp system with
  singular potential},
\newblock {Z. Anal. Anwend.}, \textbf{25}, (2006), 51--72.

\bibitem{MiranvilleQuintanilla09} (MR2524435)
A.~Miranville and R.~Quintanilla,
\newblock \emph{A generalization of the {C}aginalp phase-field system based on the
  {C}attaneo law},
\newblock {Nonlinear Anal.}, \textbf{71}, (2009), 2278--2290.

\bibitem{MiranvilleQuintanilla10} (MR2661949)
A.~Miranville and R.~Quintanilla,
\newblock \emph{A {C}aginalp phase-field system with a nonlinear coupling},
\newblock {Nonlinear Anal. Real World Appl.}, \textbf{11}, (2010), 2849--2861.

\bibitem{Ruyer}
Pierre Ruyer,
\newblock ``Mod\`ele de champ de phase pour {l'}\'etude de
  {l'}\'ebullition",
\newblock Ph.D thesis, \'Ecole Polytechnique, 2006.

\end{thebibliography}
\end{document}